\documentclass[final]{siamart1116}
\usepackage{amssymb}
\usepackage{eurosym}
\usepackage{lipsum}
\usepackage{amsfonts}
\usepackage{graphicx}
\usepackage{epstopdf}
\usepackage{algorithmic}
\usepackage{amsopn}

\ifpdf
  \DeclareGraphicsExtensions{.eps,.pdf,.png,.jpg}
\else
  \DeclareGraphicsExtensions{.eps}
\fi
\numberwithin{theorem}{section}
\newcommand{\TheTitle}{Convexification of a Parabolic Inverse
Problem}
\newcommand{\TheAuthors}{M.V.Klibanov J.Li and W.Zhang}

\headers{\TheTitle}{\TheAuthors}

\ifpdf
\hypersetup{
  pdftitle={\TheTitle},
  pdfauthor={\TheAuthors}
}
\fi

\renewcommand{\thefigure}{\arabic{section}.\arabic{figure}}
\renewcommand\theequation{\thesection.\arabic{equation}}

\begin{document}

\title{Convexification for an Inverse Parabolic Problem\thanks{%
Submitted to the editors DATE. 
\funding{The work of Klibanov was supported by US Army Research Laboratory and US Army Research Office grant W911NF-19-1-0044.
The work of Li  was partially supported by the NSF of China No. 11971221 and 11731006, 
and Guangdong Provincial Key Laboratory of Computational Science
and Material Design No. 2019B030301001. The work of Zhang was partially supported by  the Shenzhen Sci-Tech Fund No. JCYJ20170818153840322, JCYJ20180307151603959 and  the NSF of China No. 11901282}%
}}
\author{ Michael V. Klibanov\thanks{%
Department of Mathematics and Statistics, University of North Carolina at
Charlotte, Charlotte, NC 28223, USA (mklibanv@uncc.edu)} \and Jingzhi Li%
\thanks{%
Department of Mathematics, SUSTech International Center for Mathematics and
Guangdong Provincial Key Laboratory of Computational Science and Material
Design, Southern University of Science and Technology (SUSTech), Shenzhen,
Guangdong Province, P.R.China (li.jz@sustech.edu.cn)} \and Wenlong Zhang%
\thanks{%
Department of Mathematics, Southern University of Science and Technology
(SUSTech), Shenzhen, Guangdong Province, P.R.China (zhangwl@sustech.edu.cn)} 
}
\date{}
\maketitle

\begin{abstract}
A convexification-based numerical method for a Coefficient Inverse Problem
for a parabolic PDE is presented. The key element of this method is the
presence of the so-called Carleman Weight Function in the numerical scheme.
Convergence analysis ensures the global convergence of this method, as
opposed to the local convergence of the conventional least squares
minimization techniques. Numerical results demonstrate a good performance.
\end{abstract}

\makeatletter
\@addtoreset{equation}{section} \makeatother



\section{Introduction}

\label{sec:1}

In this paper, we construct a globally convergent numerical method for a
Coefficient Inverse Problem (CIP) for a parabolic PDE. This method is \
based on the so-called \emph{convexification} concept. Both convergence
analysis and numerical results are presented. The CIP, which is considered
here, has applications in heat conduction \cite{Alifanov}, diffusion theory 
\cite{Locarx} and in medical optical imaging using the diffuse infrared
light \cite{Das}. In addition, this CIP has applications in financial
mathematics in the search of the volatility coefficient in the Black-Scholes
equation using the market data \cite{Is1,KBS}. In the latter case, the
volatility coefficient should be assumed to be dependent on the stock price.

The most challenging question one needs to address prior a numerical
treatment of any CIP for a PDE is: \emph{How to choose such a starting point
of iterations that the convergence of a corresponding iterative numerical
method to the correct solution of that CIP would be \textbf{rigorously
guaranteed}?} The underlying reason of the importance of this question is
that CIPs for PDEs are both nonlinear and ill-posed. These two factors cause
the well known phenomenon of multiple local minima and ravines of
conventional Tikhonov-like least squares cost functionals for CIPs, see,
e.g. \cite{Scales} for a convincing numerical example of this phenomenon.
Therefore, the above question cannot be addressed within the framework of
the conventional least squares minimization.

One option would be to choose that starting point in a small neighborhood of
the solution. However, such a good first guess is rarely available in
applications. In fact, in such a case, the rest of the numerical procedure
would be a \emph{locally} convergent numerical method. On the other hand, we
call a numerical method for a CIP \emph{globally convergent}, if there
exists a theorem claiming that this method delivers at least one point in a
sufficiently small neighborhood of the exact solution without any advanced
knowledge of this neighborhood.

To address the above question, the first author with coauthors has been
working since 1995 \cite{BK,Klib95,Klib97,KT,KK} on the concept of the
so-called \emph{convexification} method for CIPs. This concept leads to
globally convergent numerical methods. Those initial works on the
convexification were not concerned with numerical studies (although, see 
\cite{KT} for some numerical results in the 1D case). The main reason of
this was the lack of some theorems at that time, which would ensure a proper
behavior of iterates. These theorems were first proved in \cite{Bak}.

After \cite{Bak}, a number of works on the convexification was published by
the first author with coauthors, in which the theory is combined with
numerical results, see, e.g. \cite{Khoa,convIPnew,EIT,timedomain,KlibKol3}.
We also refer to \cite{Baud} where a different version of the
convexification is developed for a CIP for the hyperbolic equation $%
u_{tt}=\Delta u+q\left( x\right) u$ and numerical results are presented.
Most recently the idea of \cite{Baud} was explored in \cite{Boulakia,LeLoc}
to develop globally convergent numerical methods for some inverse problems
for quasilinear parabolic PDEs. We also refer to the most recent work \cite%
{Harrach} for another idea of a globally convergent numerical method for a
discrete statement of a special version of the electrical impedance
tomography problem.

The convexification is a concept rather than a ready-to-use algorithm. This
means that each new CIP\ requires it own version of the convexification, and
these versions differ from each other quite significantly. Currently the
convexification is developed analytically and tested numerically for CIPs
for the Helmholtz equation \cite{Khoa,convIPnew,KlibKol3}, two hyperbolic
equations \cite{Baud,BK,timedomain} and Electrical Impedance Tomography \cite%
{EIT}. The goal of this paper is to develop analytically and implement
computationally the convexification method for a CIP for a parabolic PDE.
The first step towards this goal was made in \cite{KK}. However, there are
some problems in \cite{KK}, which prevent one from a numerical
implementation of the idea of \cite{KK}. Indeed, although a weighted
globally strictly convex Tikhonov-like functional is constructed in \cite{KK}%
, the Carleman Weight Function (CWF) in it is too complicated since it
depends on two large parameters rather than on a single one. This means that
the CWF of \cite{KK} changes too rapidly. The latter does not allow a
numerical implementation, see \cite{Baud} for a similar conclusion regarding
a different CIP. In addition, since \cite{KK} was published before \cite{Bak}%
, then uniqueness and existence of the minimizer as well as the global
convergence of the gradient projection method are not proven in \cite{KK}.
Besides, numerical studies were not conducted in \cite{KK}.

Thus, in this paper we first prove a new Carleman estimate with a simpler
CWF, which can be used for computations. Next, we prove the central result:
the global strict convexity of our weighted Tikhonov-like functional. Next,
we establish the existence and uniqueness of its minimizer, estimate the
distance between that minimizer and the exact solution and prove the global
convergence of the gradient projection method to the exact solution.
Finally, we describe results of our numerical experiments.

In the convexification, one constructs a weighted Tikhonov-like functional $%
J_{\lambda }$, where $\lambda \geq 1$ is the parameter. The weight is the
CWF, i.e. the function which is involved as the weight in the Carleman
estimate for the underlying PDE operator. Given a convex bounded set $%
B\left( d\right) \subset H^{k}$ of an arbitrary diameter $d>0$ in a certain
Hilbert space $H^{k}$, one can choose the parameter $\lambda $ of the CWF
such that the strict convexity of that functional on $\overline{B\left(
d\right) }$ is ensured. Thus, the local minima do not exist. Furthermore, as
stated above, starting from the publication \cite{Bak}, all works about the
convexification contain theorems, which claim the existence and uniqueness
of the minimizer of $J_{\lambda }$ on the set $\overline{B\left( d\right) }$
and convergence of the gradient projection method of the minimization of $%
J_{\lambda }$ to that minimizer, if starting from an \emph{arbitrary} point
of $B\left( d\right) .$ Next, as long as the level of the noise in the data
tends to zero, those minimizers converge to the correct solution of the
corresponding CIP. In particular, the latter means the stability of
minimizers with respect to a small noise in the data. Since the diameter $d$
of the convex set $B\left( d\right) $ is an arbitrary one, then the latter
amounts to the \emph{global convergence}. Even though the theory requires
the parameter $\lambda $ to be sufficiently large, our rich computational
experience with the convexification shows that in real computations $\lambda
\in \left[ 1,3\right] $ is sufficient \cite%
{Khoa,convIPnew,EIT,timedomain,KlibKol3}, also, see (\ref{9.1}). In other
words, computations are far less pessimistic than the theory is.

In section 2, we formulate both forward and inverse problems. The first step
of section 3 consists in obtaining a nonlinear integral differential
equation in which the unknown coefficient is not present. In the second step
of that section we construct the above mentioned weighted Tikhonov-like
functional with a CWF in it. In section 4 we formulate our theorems related
to this functional. These theorems are proved in sections 5-8. In section 9,
we present results of our numerical studies.

\section{Statement of the Coefficient Inverse Problem}

\label{sec:2}

Below $\mathbf{x}=\left( x,\overline{x}\right) \in \mathbb{R}^{n},$ where $%
\overline{x}=\left( x_{2},...,x_{n}\right) $ and $x=x_{1}.$ Let the numbers $%
A,B>0$ and $A<B.$ We introduce the cube $\Omega \subset \mathbb{R}^{n}$ and
a part $\Gamma $ of its boundary $\partial \Omega $ as%
\begin{equation}
\Omega =\left\{ \mathbf{x}:A<x,x_{2},...,x_{n}<B\right\} ,\Gamma =\left\{
x=B,A<x_{2},...,x_{n}<B\right\} .  \label{2.1}
\end{equation}%
Let the number $T>0.$ Denote%
\[
Q_{T}^{\pm }=\Omega \times \left( -T,T\right) ,S_{T}^{\pm }=\partial \Omega
\times \left( -T,T\right) ,\Gamma _{T}^{\pm }=\Gamma \times \left(
-T,T\right) . 
\]%
Below $\alpha \in \left( 0,1\right) ,m\geq 1$ is an integer and $C^{m+\alpha
}\left( \overline{\Omega }\right) ,C^{2m+\alpha ,m+\alpha /2}\left( 
\overline{Q_{T}^{\pm }}\right) $ are H\"{o}lder spaces \cite{Lad}. Let%
\[
b_{j}\left( \mathbf{x}\right) ,c\left( \mathbf{x}\right) \in C^{2+\alpha
}\left( \overline{\Omega }\right) ;\text{ }j=1,...,n. 
\]%
We consider the elliptic operator $L$ in the following form: 
\begin{equation}
Lu=\Delta u+\mathop{\displaystyle \sum }\limits_{j=1}^{n}b_{j}\left( \mathbf{%
x}\right) u_{x_{j}}-c\left( \mathbf{x}\right) u,\text{ }\mathbf{x}\in \Omega
.  \label{2.2}
\end{equation}%
We assume that 
\begin{equation}
c\left( \mathbf{x}\right) \geq 0\text{ \ in }\overline{\Omega }.
\label{2.30}
\end{equation}

The forward parabolic initial boundary value problem is stated as \cite{Lad}:

\textbf{Forward Problem}. \emph{Let the initial condition} $f\left( \mathbf{x%
}\right) \in C^{4+\alpha }\left( \overline{\Omega }\right) .$\emph{\ Find a
function} $u\left( \mathbf{x},t\right) \in C^{4+\alpha ,2+\alpha /2}\left( 
\overline{Q_{T}^{\pm }}\right) $ \emph{satisfying the following conditions:}%
\begin{equation}
u_{t}=Lu\text{ in }Q_{T}^{\pm },  \label{2.4}
\end{equation}%
\begin{equation}
u\left( \mathbf{x},-T\right) =f\left( \mathbf{x}\right) ,  \label{2.5}
\end{equation}%
\begin{equation}
u\mid _{S_{T}^{\pm }}=g_{0}\left( \mathbf{x},t\right) .  \label{2.6}
\end{equation}

If the domain $\Omega $ would have its boundary $\partial \Omega \in
C^{4+\alpha }$ and if the Dirichlet condition $g_{0}\left( \mathbf{x}%
,t\right) $ would belong to $C^{4+\alpha ,2+\alpha /2}\left( \overline{%
S_{T}^{\pm }}\right) $ and also corresponding compatibility conditions would
be satisfied \cite{Lad}, then the existence and uniqueness of the solution $%
u\in C^{4+\alpha ,2+\alpha /2}\left( \overline{Q_{T}^{\pm }}\right) $ of
problem (\ref{2.2})-(\ref{2.6}) would be ensured \cite{Lad}. However, for
the the convenience of our derivations for the inverse problem, we have
chosen the case of a piecewise smooth boundary $\partial \Omega .$ Hence, we
can only assume the existence of the solution $u\in C^{4+\alpha ,2+\alpha
/2}\left( \overline{Q_{T}^{\pm }}\right) $ of problem (\ref{2.4})-(\ref{2.6}%
).\ As to its uniqueness, it follows immediately from (\ref{2.30}) and the
maximum principle for parabolic PDEs.

\textbf{Coefficient Inverse Problem (CIP). }\emph{Let the number }$t_{0}\in
\left( -T,T\right) .$\emph{\ Suppose that the following two functions }$%
g_{1}\left( \mathbf{x},t\right) $\emph{\ and }$f_{0}\left( \mathbf{x}\right) 
$\emph{\ are known:}%
\begin{equation}
u_{x}\mid _{\Gamma _{T}^{\pm }}=g_{1}\left( \mathbf{x},t\right) ,
\label{2.7}
\end{equation}%
\begin{equation}
u\left( \mathbf{x},t_{0}\right) =f_{0}\left( \mathbf{x}\right) .  \label{2.8}
\end{equation}%
\emph{Find the unknown coefficient }$c\left( \mathbf{x}\right) .$

If $n=3$ and functions $b_{j}\left( \mathbf{x}\right) \equiv 0$ for $%
j=1,...,n,$ then $c\left( \mathbf{x}\right) $ is the absorption coefficient
in the case of medical optical imaging using the diffuse infrared light \cite%
{Das}. Uniqueness of this CIP for any value of $T$ was proven by the first
author using the method of \cite{BukhKlib}, see, e.g. theorem 1.10.7 in \cite%
{BKbook}, theorem 2 in \cite{Klib84}, theorem 3.10 in \cite{Klib92} and
theorem 3.4 in \cite{Ksurvey}. We also refer to \cite{IY,Yam} for the
Lipschitz stability estimate for this CIP.

The data for our CIP are non redundant, so as for all CIPs for which the
convexification method works. In other words, the number $m$ of free
variables in the data equals the number $n$ of free variables in the unknown
coefficient, $m=n$. As to the globally convergent numerical methods for CIPs
with redundant data with $m>n$, see, e.g. \cite{deHoop,Kab1,Kab2}.

\section{Weighted Globally Strictly Convex Tikhonov-like Functional}

\label{sec:3}

We assume below that there exists a number $\mu >0$ such that 
\begin{equation}
f\left( \mathbf{x}\right) \geq \mu ,\text{ }\forall \mathbf{x}\in \overline{%
\Omega }\text{,}  \label{3.1}
\end{equation}%
\begin{equation}
g_{0}\left( \mathbf{x},t\right) \geq \mu ,\text{ }\forall \left( \mathbf{x}%
,t\right) \in \overline{S_{T}^{\pm }}.  \label{3.100}
\end{equation}%
Then (\ref{2.30}), (\ref{3.1}), (\ref{3.100}) and the maximum principle for
parabolic PDEs \cite{Lad} imply that 
\begin{equation}
u\left( \mathbf{x},t\right) \mathbf{\geq }\mu \text{ in }\overline{%
Q_{T}^{\pm }}\mathbf{.}  \label{3.2}
\end{equation}

\subsection{Nonlinear integral differential equation}

\label{sec:3.1}

Using (\ref{3.2}), we introduce a new function $v\left( \mathbf{x},t\right)
, $ 
\begin{equation}
v\left( \mathbf{x},t\right) =\ln u\left( \mathbf{x},t\right) \rightarrow
u=e^{v}.  \label{3.3}
\end{equation}%
Substituting (\ref{3.3}) in (\ref{2.4})-(\ref{2.8}), we obtain in $%
Q_{T}^{\pm }:$%
\begin{equation}
v_{t}-\Delta v-\left( \nabla v\right) ^{2}-\mathop{\displaystyle \sum }%
\limits_{k=1}^{n}b_{j}\left( \mathbf{x}\right) v_{x_{j}}=c\left( \mathbf{x}%
\right) ,  \label{3.4}
\end{equation}%
\begin{equation}
v\mid _{S_{T}^{\pm }}=\ln g_{0}\left( \mathbf{x},t\right) ,v_{x}\mid
_{\Gamma _{T}^{\pm }}=\left( g_{1}/g_{0}\right) \left( \overline{x},t\right)
,  \label{3.5}
\end{equation}%
\begin{equation}
v\left( \mathbf{x},t_{0}\right) =\ln f_{0}\left( \mathbf{x}\right) :=%
\widetilde{f}_{0}\left( \mathbf{x}\right) .  \label{3.6}
\end{equation}%
For brevity we set below $t_{0}:=0.$ The case $t_{0}\neq 0$ can be
considered along the same lines. Differentiate both sides of the nonlinear
equation (\ref{3.4}) with respect to $t$ and denote $w\left( \mathbf{x}%
,t\right) =v_{t}\left( \mathbf{x},t\right) .$ Since the function $c\left( 
\mathbf{x}\right) $ is independent on $t$, then the right hand side of the
resulting equation will be zero. By (\ref{3.6}) 
\begin{equation}
v\left( \mathbf{x},t\right) =\mathop{\displaystyle \int}\limits_{0}^{t}w%
\left( \mathbf{x},\tau \right) d\tau +\widetilde{f}_{0}\left( \mathbf{x}%
\right) ,\text{ }\left( \mathbf{x},t\right) \in Q_{T}^{\pm }.  \label{3.7}
\end{equation}%
Substituting (\ref{3.7}) in (\ref{3.4}) and (\ref{3.5}), we obtain a
nonlinear integral differential PDE with Volterra integrals, supplied by the
lateral Cauchy data,%
\[
K\left( w\right) =w_{t}-\Delta w-\mathop{\displaystyle \sum }%
\limits_{j=1}^{n}b_{j}\left( \mathbf{x}\right) w_{x_{j}} 
\]%
\begin{equation}
-2\nabla w\mathop{\displaystyle \int}\limits_{0}^{t}\nabla w\left( \mathbf{x}%
,\tau \right) d\tau -2\nabla w\nabla \widetilde{f}_{0}=0,\text{ }\left( 
\mathbf{x},t\right) \in Q_{T}^{\pm },  \label{3.8}
\end{equation}%
\begin{equation}
w\mid _{S_{T}^{\pm }}=p_{0}\left( \mathbf{x},t\right) ,w_{x}\mid _{\Gamma
_{T}^{\pm }}=p_{1}\left( \mathbf{x},t\right) ,  \label{3.9}
\end{equation}%
where $p_{0}\left( \mathbf{x},t\right) =\left( g_{0t}/g_{0}\right) \left( 
\mathbf{x},t\right) $ and $p_{1}\left( \mathbf{x},t\right) =\partial
_{t}\left( g_{1}/g_{0}\right) \left( \overline{x},t\right) .$

\subsection{The functional}

\label{sec:3.2}

There are many possible choices of the CWF for the parabolic operator, see,
e.g. \cite{BKbook,IY,Ksurvey,LRS,Yam}. However, among all these choices, we
should select such a CWF\ which would be simple and would work well
computationally. Indeed, for example, the CWF of \cite{BKbook,Ksurvey,LRS}
depends on two large parameters, which means that it changes too rapidly. As
it was stated in Introduction, that rapid change prevents one from a
numerical implementation. Thus, we have chosen the CWF $\varphi _{\lambda
}\left( x,t\right) $ as:%
\begin{equation}
\varphi _{\lambda }\left( x,t\right) =\exp \left( 2\lambda \left(
x^{2}-t^{2}\right) \right) ,  \label{3.10}
\end{equation}%
where $\lambda \geq 1$ is a parameter. This means that we need to prove the
Carleman estimate with this CWF, see Theorem 1 in section 4. Let $\left[
\left( n+1\right) /2\right] $ be the maximal integer which does not exceed $%
\left( n+1\right) /2.$ Denote $k_{n}=\left[ \left( n+1\right) /2\right] +2.$
For example, we have for most popular cases of $n=1,2,3:$ 
\[
k_{n}=\left\{ 
\begin{array}{c}
3\text{ if }n=1,2, \\ 
4\text{ if }n=3\text{.}%
\end{array}%
\right. 
\]%
We have chosen the number $k_{n}$ in such a way that 
\begin{equation}
H^{k_{n}}\left( Q_{T}^{\pm }\right) \subseteq H^{3}\left( Q_{T}^{\pm
}\right) ,  \label{3.11}
\end{equation}%
\begin{equation}
H^{k_{n}}\left( Q_{T}^{\pm }\right) \subset C^{1}\left( \overline{Q_{T}^{\pm
}}\right) ,\left\Vert q\right\Vert _{C^{1}\left( \overline{Q_{T}^{\pm }}%
\right) }\leq C_{0}\left\Vert q\right\Vert _{H^{k_{n}}\left( Q_{T}^{\pm
}\right) },\forall q\in H^{k_{n}}\left( Q_{T}^{\pm }\right) ,  \label{3.12}
\end{equation}%
where the number $C_{0}=C_{0}\left( Q_{T}^{\pm }\right) >0$ depends only on
the domain $Q_{T}^{\pm }.$ Relations (\ref{3.12}) follow from (\ref{3.11})
and the embedding theorem.

Let $R>0$ be an arbitrary number. We define the bounded set of functions $%
B\left( R,p_{0},p_{1}\right) $ as follows:%
\begin{eqnarray}
&&B\left( R,p_{0},p_{1}\right)  \label{3.13} \\
&=&\left\{ w\in H^{k_{n}}\left( Q_{T}^{\pm }\right) :\left\Vert w\right\Vert
_{H^{k_{n}}\left( Q_{T}^{\pm }\right) }<R,w\mid _{S_{T}^{\pm
}}=p_{0},w_{x}\mid _{\Gamma _{T}^{\pm }}=p_{1}\right\} ,  \nonumber
\end{eqnarray}%
where functions $p_{0},p_{1}$ are taken from (\ref{3.9}).

Let $\beta >0$ be a small regularization parameter and $K\left( w\right) $
be the nonlinear integral differential operator defined in (\ref{3.8}). We
construct our weighted Tikhonov-like functional with the CWF (\ref{3.10}) in
it as:%
\begin{equation}
J_{\lambda ,\beta }\left( w\right) =e^{-2\lambda B^{2}}\mathop{\displaystyle
\int}\limits_{Q_{T}^{\pm }}\left( K\left( w\right) \right) ^{2}\varphi
_{\lambda }d\mathbf{x}dt+\beta \left\Vert w\right\Vert _{H^{k_{n}}\left(
Q_{T}^{\pm }\right) }^{2}.  \label{3.15}
\end{equation}%
Since $\max_{\overline{Q_{T}^{\pm }}}\varphi _{\lambda }=e^{2\lambda B^{2}},$
then the multiplier $e^{-2\lambda B^{2}}$ is introduced in (\ref{3.15}) to
balance two terms in the right hand side of (\ref{3.15}).

\textbf{Minimization Problem}. \emph{Minimize the functional }$J_{\lambda
,\beta }\left( w\right) $\emph{\ on the set }$\overline{B\left( R\right) }$%
\emph{\ \ defined in (\ref{3.13}).}

Assume for a moment that a minimizer $w_{\min ,\lambda ,\beta }\left( 
\mathbf{x},t\right) $ of functional (\ref{3.15}) exists and is computed.
Then we first calculate the corresponding function $v_{\text{comp}}\left( 
\mathbf{x},t\right) $ via (\ref{3.7}). Next, substituting $v_{\min ,\lambda
,\beta }\left( \mathbf{x},t\right) =$ in equation (\ref{3.4}), we calculate
an approximation for the target unknown coefficient $c\left( \mathbf{x}%
\right) .$ However, due to the inevitable computational errors as well as
the noise in the data, the resulting left hand side of (\ref{3.4}) would
depend on $t$. Hence, to calculate an approximation $c_{\text{comp}}\left( 
\mathbf{x}\right) $ for $c\left( \mathbf{x}\right) ,$ we set%
\begin{equation}
c_{\text{comp}}\left( \mathbf{x}\right) =  \label{3.16}
\end{equation}%
\[
\frac{1}{2\gamma T}\mathop{\displaystyle \int}\limits_{-\gamma T}^{\gamma
T}\left( \partial _{t}v_{\text{comp}}-\Delta v_{\text{comp}}-\left( \nabla
v_{\text{comp}}\right) ^{2}-\mathop{\displaystyle \sum }%
\limits_{k=1}^{n}b_{j}\left( \mathbf{x}\right) \partial _{x_{j}}v_{\text{comp%
}}\right) dt, 
\]%
where the number $\gamma \in \left( 0,1/\sqrt{3}\right) $ is chosen in
section 4. Thus, we focus below on the Minimization Problem.

\section{Theorems}

\label{sec:4}

Introduce the subspaces $H_{0}^{2,1}\left( Q_{T}^{\pm }\right) \subset $ $%
H^{2,1}\left( Q_{T}^{\pm }\right) $ and $H_{0}^{k_{n}}\left( Q_{T}^{\pm
}\right) \subset H^{k_{n}}\left( Q_{T}^{\pm }\right) $ as%
\[
H_{0}^{2,1}\left( Q_{T}^{\pm }\right) =\left\{ u\in H^{2,1}\left( Q_{T}^{\pm
}\right) :u\mid _{S_{T}^{\pm }}=0,u_{x}\mid _{\Gamma _{T}^{\pm }}=0\right\}
, 
\]%
\[
H_{0}^{k_{n}}\left( Q_{T}^{\pm }\right) =\left\{ u\in H^{k_{n}}\left(
Q_{T}^{\pm }\right) :u\mid _{S_{T}^{\pm }}=0,u_{x}\mid _{\Gamma _{T}^{\pm
}}=0\right\} . 
\]

Since it is well known that any Carleman estimate depends only on the
principal part of the operator, see, e.g. \cite{Ksurvey,LRS}, then we
consider in Theorem 1 only the principal part $\partial _{t}-\Delta $ of the
parabolic operator $\partial _{t}-L.$

\textbf{Theorem 1} (Carleman estimate).\emph{\ Suppose that the domain }$%
\Omega $\emph{\ and the CWF }$\varphi _{\lambda }\left( x,t\right) $\emph{\
are the same as in (\ref{2.1}) and (\ref{3.10}) respectively. Then there
exist numbers }$\lambda _{0},C,$\emph{\ }%
\begin{equation}
\lambda _{0}=\lambda _{0}\left( \Omega \right) \geq 1\emph{,}\text{ }%
C=C\left( \Omega ,T\right) >0\emph{\ }  \label{4.02}
\end{equation}%
\emph{depending only on listed parameters such that the following Carleman
estimate holds}%
\[
\mathop{\displaystyle \int}\limits_{Q_{T}^{\pm }}\left( u_{t}-\Delta
u\right) ^{2}\varphi _{\lambda }d\mathbf{x}dt\geq \frac{C}{\lambda }%
\mathop{\displaystyle \int}\limits_{Q_{T}^{\pm }}\left( u_{t}^{2}+%
\mathop{\displaystyle \sum }\limits_{i,j=1}^{n}u_{x_{i}x_{j}}^{2}\right)
\varphi _{\lambda }d\mathbf{x}dt 
\]%
\begin{equation}
+C\lambda \mathop{\displaystyle \int}\limits_{Q_{T}^{\pm }}\left[ \left(
\nabla u\right) ^{2}+\lambda ^{2}u^{2}\right] \varphi _{\lambda }d\mathbf{x}%
dt  \label{4.2}
\end{equation}%
\[
-C\exp \left( 2\lambda \left( B^{2}-T^{2}\right) \right) \mathop{%
\displaystyle \int}\limits_{\Omega }\left[ \left( u_{t}^{2}+\left( \nabla
u\right) ^{2}+\lambda ^{2}u^{2}\right) \left( \mathbf{x},T\right) \right] d%
\mathbf{x} 
\]%
\[
-C\exp \left( 2\lambda \left( B^{2}-T^{2}\right) \right) \mathop{%
\displaystyle \int}\limits_{\Omega }\left[ \left( u_{t}^{2}+\left( \nabla
u\right) ^{2}+\lambda ^{2}u^{2}\right) \left( \mathbf{x},-T\right) \right] d%
\mathbf{x,} 
\]%
\[
\forall \lambda \geq \lambda _{0},\forall u\in H_{0}^{2,1}\left( Q_{T}^{\pm
}\right) . 
\]

\textbf{Remarks 1:}

\textbf{1.} \emph{An analog of estimate (\ref{4.2}) was proven in \cite{KKLY}%
, although only for the 1D case, and terms with }$u_{xx}$\emph{,}$u_{t}$%
\emph{\ were not involved in the estimate of \cite{KKLY}. However, the
presence in (\ref{4.2}) of the terms with derivatives involved in the
principal part of the parabolic operator is important for the proofs of
Theorems 2-6. Thus, Carleman estimate (\ref{4.2}) is new.}

\textbf{2.} \emph{Since the normal derivative of the function }$u\in
H_{0}^{2,1}\left( Q_{T}^{\pm }\right) $\emph{\ equals zero only on the part }%
$\Gamma _{T}^{\pm }$\emph{\ of the lateral boundary }$S_{T}^{\pm }$\emph{\
of the time cylinder }$Q_{T}^{\pm }$\emph{\ rather than on the whole }$%
S_{T}^{\pm },$\emph{\ then one should carefully analyze integrals over }$%
S_{T}^{\pm }$\emph{\ which occur in the pointwise Carleman estimate: to make
sure that these integrals equal zero.}

\textbf{Theorem 2} (the central theorem of this paper). \emph{Assume that
condition (\ref{3.2}) holds. The functional }$J_{\lambda ,\beta }\left(
w\right) $\emph{\ has the Frech\'{e}t derivative }$J_{\lambda ,\beta
}^{\prime }\left( w\right) \in H_{0}^{k_{n}}\left( Q_{T}^{\pm }\right) $ 
\emph{\ for all }$\lambda ,\beta >0,w\in B\left( 3R,p_{0},p_{1}\right) .$%
\emph{\ Let }$\lambda _{0}\geq 1$\emph{\ be the constant of Theorem 1. There
exist constants }%
\begin{equation}
\lambda _{1}=\lambda _{1}\left( R,A,B,T,\max_{j}\left\Vert b_{j}\right\Vert
_{C\left( \overline{\Omega }\right) },\left\Vert f_{0}\right\Vert
_{C^{1}\left( \overline{\Omega }\right) },\mu \right) \geq \lambda _{0},
\label{4.200}
\end{equation}%
\begin{equation}
C_{1}=C_{1}\left( R,A,B,T,\max_{j}\left\Vert b_{j}\right\Vert _{C\left( 
\overline{\Omega }\right) },\left\Vert f_{0}\right\Vert _{C^{1}\left( 
\overline{\Omega }\right) },\mu \right) >0  \label{4.20}
\end{equation}%
\emph{\ depending only on listed parameters such that if }$\lambda \geq
\lambda _{1}$ \emph{and the regularization parameter }$\beta \in \left[
2e^{-\lambda T^{2}},1\right) ,$\emph{\ then the functional }$J_{\lambda
,\beta }\left( w\right) $\emph{\ is strictly convex on the set }$\overline{%
B\left( R,p_{0},p_{1}\right) }$ \emph{\ for all }$\lambda \geq \lambda _{1},$%
\emph{\ i.e. for all} $w_{1},w_{2}\in \overline{B\left( R,p_{0},p_{1}\right) 
}$ \emph{and for all }$\lambda \geq \lambda _{1}$ 
\begin{equation}
J_{\lambda ,\beta }\left( w_{2}\right) -J_{\lambda ,\beta }\left(
w_{1}\right) -J_{\lambda ,\beta }^{\prime }\left( w_{1}\right) \left(
w_{2}-w_{1}\right)  \label{4.201}
\end{equation}%
\[
\geq \frac{C_{1}}{\lambda }\exp \left( -2\lambda \left(
T^{2}+B^{2}-A^{2}\right) \right) \left\Vert w_{2}-w_{1}\right\Vert
_{H^{2,1}\left( Q_{T}^{\pm }\right) }^{2}+\frac{\beta }{2}\left\Vert
w_{2}-w_{1}\right\Vert _{H^{k_{n}}\left( Q_{T}^{\pm }\right) }^{2}. 
\]

Everywhere below $C>0$ and $C_{1}>0$ denote different constants depending
only on parameters listed in (\ref{4.02}) and (\ref{4.20}) respectively.

\textbf{Theorem 3}. \emph{Assume that condition (\ref{3.2}) holds. Let
parameters }$\lambda _{1},\lambda \geq \lambda _{1}$ \emph{and }$\beta $%
\emph{\ be the same as the ones in Theorem 2. Then there exists unique
minimizer }$w_{\min ,\lambda ,\beta }\in \overline{B\left( R\right) }$\emph{%
\ of the functional }$J_{\lambda ,\beta }\left( w\right) $\emph{\ on the set 
}$\overline{B\left( R\right) }.$\emph{\ Furthermore, the following
inequality holds:}%
\[
J_{\lambda ,\beta }^{\prime }\left( w_{\min ,\lambda ,\beta }\right) \left(
w-w_{\min ,\lambda ,\beta }\right) \geq 0,\text{ }\forall w\in \overline{%
B\left( R\right) }. 
\]

Following the regularization theory \cite{T}, we assume now that there
exists an ideal, the so-called `exact' solution $c^{\ast }\left( \mathbf{x}%
\right) \in C^{2+\alpha }\left( \overline{\Omega }\right) $ of the CIP (\ref%
{2.30}), (\ref{2.4})-(\ref{2.8}), where the data (\ref{2.6})-(\ref{2.8}) are
noiseless. Also, let $c_{\text{comp}}\left( \mathbf{x}\right) $ be the
coefficient $c\left( \mathbf{x}\right) $ reconstructed from the minimizer $%
w_{\min ,\lambda ,\beta }\left( \mathbf{x},t\right) $ via backwards
calculations, as outlined in the last paragraph of section 3 and, in
particular, in (\ref{3.16}). Having the function $c^{\ast }\left( \mathbf{x}%
\right) ,$ one can construct the noise free solution $w^{\ast }\in
H^{k_{n}}\left( Q_{T}^{\pm }\right) $ of equation (\ref{3.8}) with the
noiseless boundary data $p_{0}^{\ast },p_{1}^{\ast }$ in (\ref{3.9}) and the
noiseless function $\widetilde{f}_{0}^{\ast }\left( \mathbf{x}\right) $ in (%
\ref{3.8}).

We now want to estimate the distance between the minimizer $w_{\min ,\lambda
,\beta }$ and the function $w^{\ast }$ as well as between coefficients $%
c^{\ast }\left( \mathbf{x}\right) $ and $c_{\text{comp}}\left( \mathbf{x}%
\right) .$ To do this, we first arrange zero boundary conditions in an
analog of (\ref{3.9}). More precisely, we assume that there exist functions $%
G\left( \mathbf{x},t\right) $ and $G^{\ast }\left( \mathbf{x},t\right) $
satisfying the same boundary conditions as those for $w$ and $w^{\ast }$
respectively and such that their norms in the space $H^{k_{n}}\left(
Q_{T}^{\pm }\right) $ are less than $R$, i.e. 
\begin{equation}
G\in B\left( R,p_{0},p_{1}\right) ,G^{\ast }\in B\left( R,p_{0}^{\ast
},p_{1}^{\ast }\right)  \label{4.3}
\end{equation}%
Let a small number $\delta \in \left( 0,1\right) $ be the level of the noise
in the functions $G$ and $f_{0}.$ More, precisely, we assume that%
\begin{equation}
\left\Vert G-G^{\ast }\right\Vert _{H^{k_{n}}\left( Q_{T}^{\pm }\right)
}<\delta ,  \label{4.5}
\end{equation}%
\begin{equation}
\left\Vert f_{0}-f_{0}^{\ast }\right\Vert _{C^{1}\left( \overline{\Omega }%
\right) }<\delta .  \label{4.50}
\end{equation}

\textbf{Remark 2}.\emph{\ By (\ref{4.50}), we replace below }$\left\Vert
f_{0}\right\Vert _{C^{1}\left( \overline{\Omega }\right) }$\emph{\ with }$%
\left\Vert f_{0}^{\ast }\right\Vert _{C^{1}\left( \overline{\Omega }\right)
} $\emph{in (\ref{4.200}) and (\ref{4.20}).}

We also assume that functions 
\begin{equation}
w^{\ast }\in B\left( R-\delta ,p_{0}^{\ast },p_{1}^{\ast }\right) ,
\label{4.6}
\end{equation}%
\begin{equation}
\left\Vert f_{0}^{\ast }\right\Vert _{C^{1}\left( \overline{\Omega }\right)
}<R,\text{ \ }\min_{\overline{\Omega }}f_{0}^{\ast }\geq \mu >0,
\label{4.60}
\end{equation}%
where the number $\mu $ is the same as in (\ref{3.1}), (\ref{3.2}) and is
independent on $\delta $. Then (\ref{3.9}), (\ref{4.5}) and (\ref{4.6})
imply that 
\begin{equation}
G\in B\left( R,p_{0},p_{1}\right) .  \label{4.7}
\end{equation}%
Denote 
\begin{equation}
W=w-G,W^{\ast }=w^{\ast }-G^{\ast }.  \label{4.8}
\end{equation}%
Similarly with (\ref{3.13}) denote%
\begin{equation}
B_{0}\left( 2R\right) =\left\{ W\in H^{k_{n}}\left( Q_{T}^{\pm }\right)
:\left\Vert W\right\Vert _{H^{k_{n}}\left( Q_{T}^{\pm }\right) }<2R,W\mid
_{S_{T}^{\pm }}=W_{x}\mid _{\Gamma _{T}^{\pm }}=0\right\} .  \label{4.13}
\end{equation}%
Then (\ref{4.6})-(\ref{4.8}) imply that 
\begin{equation}
W\in B_{0}\left( 2R\right) ,\forall w\in B\left( R,p_{0},p_{1}\right) \text{
and also }W^{\ast }\in B_{0}\left( 2R-2\delta \right) ,  \label{4.14}
\end{equation}%
\begin{equation}
W+G\in B\left( 3R,p_{0},p_{1}\right) ,\forall W\in B_{0}\left( 2R\right) .
\label{4.15}
\end{equation}%
Due to (\ref{4.15}), it is convenient to denote below $\lambda _{1}\left(
3R\right) ,\lambda \left( 3R\right) ,$ which means that the values of the
parameters $\lambda _{1}$ and $\lambda \geq \lambda _{1}$ correspond to $%
B\left( 3R,p_{0},p_{1}\right) $ in Theorem 2 and, in particular, $R$ is
replaced with $3R$ in (\ref{4.200}) and (\ref{4.20}). Consider the
functional $I_{\lambda ,\beta }\left( W\right) ,$%
\begin{equation}
I_{\lambda ,\beta }:B_{0}\left( 2R\right) \rightarrow \mathbb{R},\text{ }%
I_{\lambda ,\beta }\left( W\right) =J_{\lambda ,\beta }\left( W+G\right) .
\label{4.16}
\end{equation}

\textbf{Theorem 4}.\ \emph{Assume that condition (\ref{3.2}) holds. Let
parameters }$\lambda _{1}$\emph{\ and }$\beta $\emph{\ be the same as in
Theorem 2, except that }$R$\emph{\ is replaced with }$3R$\emph{\ in (\ref%
{4.200}). Then the functional }$I_{\lambda ,\beta }\left( W\right) $\emph{\
is strictly convex on the ball }$\overline{B_{0}\left( 2R\right) }$ \emph{\
for all }$\lambda \geq \lambda _{1}\left( 3R\right) .$ \emph{Here,} $\lambda
_{1}\left( 3R\right) $ \emph{means (\ref{4.200}), where }$R$\emph{\ is
replaced with }$3R$\emph{\ and }$f_{0\text{ }}$\emph{is replaced with }$f_{0%
\text{ }}^{\ast }$\emph{\ (Remark 2). In other words, the following analog
of (\ref{4.201}) holds}%
\begin{equation}
I_{\lambda ,\beta }\left( W_{2}\right) -I_{\lambda ,\beta }\left(
W_{1}\right) -I_{\lambda ,\beta }^{\prime }\left( W_{1}\right) \left(
W_{2}-W_{1}\right)  \label{4.17}
\end{equation}%
\[
\geq \frac{C_{1}}{\lambda }\exp \left( -2\lambda \left(
T^{2}+B^{2}-A^{2}\right) \right) \left\Vert W_{2}-W_{1}\right\Vert
_{H^{2,1}\left( Q_{T}^{\pm }\right) }^{2}+\frac{\beta }{2}\left\Vert
W_{2}-W_{1}\right\Vert _{H^{k_{n}}\left( Q_{T}^{\pm }\right) }^{2}, 
\]%
\emph{for all }$\lambda \geq \lambda _{1}\left( 3R\right) $\emph{\ and for
all }$W_{1},W_{2}\in \overline{B_{0}\left( 2R\right) },$ \emph{where }$%
I_{\lambda ,\beta }^{\prime }\left( W\right) \in H_{0}^{k_{n}}\left(
Q_{T}^{\pm }\right) $\emph{\ is the Frech\'{e}t derivative \ of the
functional }$I_{\lambda ,\beta }\left( W\right) $\emph{\ at the point }$W$%
\emph{, which exists due to Theorem 2 and (\ref{4.16}).} \emph{Furthermore,
there exists unique minimizer }$W_{\min ,\lambda \left( 3R\right) ,\beta
}\in \overline{B_{0}\left( 2R\right) }$\emph{\ of the functional }$%
I_{\lambda ,\beta }\left( W\right) $\emph{\ and the following inequality
holds:}%
\begin{equation}
I_{\lambda \left( 3R\right) ,\beta }^{\prime }\left( W_{\min ,\lambda \left(
3R\right) ,\beta }\right) \left( W-W_{\min ,\lambda \left( 3R\right) ,\beta
}\right) \geq 0,\text{ }\forall W\in \overline{B_{0}\left( 2R\right) }.
\label{4.170}
\end{equation}

\textbf{Theorem 5} (accuracy estimates). \emph{Assume that condition (\ref%
{3.2}) holds. Suppose that conditions (\ref{4.3})-(\ref{4.8}) hold and also
let }$T>\sqrt{3\left( B^{2}-A^{2}\right) }$\emph{. Choose a number }$\gamma
\in \left( 0,1/\sqrt{3}\right) $\emph{\ such that }$T^{2}\left( 1-3\gamma
^{2}\right) >3\left( B^{2}-A^{2}\right) .$\emph{\ Denote }%
\[
\eta _{1}=\gamma ^{2}T^{2}+B^{2}-A^{2},\eta _{2}=\left( 1-3\gamma
^{2}\right) T^{2}-3\left( B^{2}-A^{2}\right) ,\rho =\frac{1}{2}\min \left( 1,%
\frac{\eta _{2}}{\eta _{1}}\right) . 
\]%
\emph{Let }$\lambda _{1}=\lambda _{1}\left( 3R\right) $\emph{\ be the number
of Theorem 4. Choose a sufficiently small number }$\delta _{0}>0$\emph{\
such that }$\ln \left( \delta _{0}^{-1/\eta _{1}}\right) \geq \lambda _{1}.$%
\emph{\ For each }$\delta \in \left( 0,\delta _{0}\right) ,$\emph{\ let }$%
\lambda =\lambda \left( \delta ,3R\right) =\ln \left( \delta ^{-1/\eta
_{1}}\right) >\lambda _{1}\left( 3R\right) .$\emph{\ Let the regularization
parameter }$\beta =\beta \left( \delta ,3R\right) =2e^{-\lambda \left(
\delta ,3R\right) T^{2}}$\emph{\ (see Theorem 2). Let }$w_{\min ,\lambda
\left( \delta ,3R\right) ,\beta \left( \delta ,3R\right) }=W_{\min ,\lambda
\left( \delta ,3R\right) ,\beta \left( \delta ,3R\right) }+G$ \emph{(Theorem
4) and let }$c_{\text{comp}}\left( \mathbf{x}\right) $\emph{\ be the
function }$c\left( \mathbf{x}\right) $ \emph{computed from the function }$%
w_{\min ,\lambda \left( \delta ,3R\right) ,\beta \left( \delta ,3R\right)
}\left( \mathbf{x},t\right) $ \emph{by the procedure described in the last
paragraph of section 3. Then the following accuracy estimates are valid}%
\begin{equation}
\left\Vert w^{\ast }-w_{\min ,\lambda \left( \delta ,3R\right) ,\beta \left(
\delta ,3R\right) }\right\Vert _{H^{2}\left( Q_{\gamma T}^{\pm }\right)
}\leq C_{2}\delta ^{\rho },  \label{4.18}
\end{equation}%
\begin{equation}
\left\Vert c^{\ast }-c_{\min ,\lambda \left( \delta ,3R\right) ,\beta \left(
\delta ,3R\right) }\right\Vert _{L_{2}\left( \Omega \right) }\leq
C_{2}\delta ^{\rho }.  \label{4.19}
\end{equation}

Here and below $C_{2}>0$ denotes different constants depending on the same
parameters as ones in (\ref{4.20}) as well as on the number $\gamma .$

We now construct the gradient projection method of the minimization of the
functional $I_{\lambda ,\beta }\left( W\right) $ defined in (\ref{4.16}) on
the set $\overline{B_{0}\left( 2R\right) }$ defined in (\ref{4.13}). Let $%
P_{B}:H^{k_{n}}\left( Q_{T}^{\pm }\right) \rightarrow \overline{B_{0}\left(
2R\right) }$ be the orthogonal projection operator. Let $W_{0}\in
B_{0}\left( 2R\right) $ be an arbitrary point of the ball $B_{0}\left(
2R\right) .$ Let the number $\omega \in \left( 0,1\right) .$ We arrange the
gradient projection method of the minimization of the functional $I_{\lambda
,\beta }\left( W\right) $ as:%
\begin{equation}
W_{n}=P_{B}\left( W_{n-1}-\omega I_{\lambda ,\beta }^{\prime }\left(
W_{n-1}\right) \right) ,\text{ }n=1,2,....  \label{4.23}
\end{equation}%
Note that since $W_{n-1},I_{\lambda ,\beta }^{\prime }\left( W_{n-1}\right)
\in H_{0}^{k_{n}}\left( Q_{T}^{\pm }\right) ,$ then the function $%
W_{n-1}-\omega I_{\lambda ,\beta }^{\prime }\left( W_{n-1}\right) $ has zero
boundary conditions (\ref{3.9}). The latter is important in the
computational practice.

\textbf{Theorem 6} (global convergence of the gradient projection method). 
\emph{Assume that condition (\ref{3.2}) holds. Let parameters }$\lambda
_{1}\left( 3R\right) $\emph{\ and }$\beta $\emph{\ be the same as in Theorem
2, except that }$R$\emph{\ is replaced with }$3R$\emph{\ in (\ref{4.200})
and let }$\lambda \left( 3R\right) $\emph{. } \emph{Then there exists a
sufficiently small number }$\omega _{0}=\omega _{0}\left( \Omega
,T,A,B,R,\lambda \right) $\emph{\ such that for any }$\omega \in \left(
0,\omega _{0}\right) $\emph{\ there exists a number }$\theta =\theta \left(
\omega \right) \in \left( 0,1\right) $\emph{\ such that the sequence (\ref%
{4.23}) converges to the unique minimizer }$W_{\min ,\lambda \left(
3R\right) ,\beta }\in \overline{B_{0}\left( 2R\right) }$\emph{\ (Theorem 4)
in the norm of the space }$H^{k_{n}}\left( Q_{T}^{\pm }\right) .$\emph{\
More precisely,}%
\begin{equation}
\left\Vert W_{\min ,\lambda \left( 3R\right) ,\beta }-W_{n}\right\Vert
_{H^{k_{n}}\left( Q_{T}^{\pm }\right) }\leq \theta ^{n}\left\Vert W_{\min
,\lambda \left( 3R\right) ,\beta }-W_{0}\right\Vert _{H^{k_{n}}\left(
Q_{T}^{\pm }\right) }.  \label{4.24}
\end{equation}

\textbf{Theorem 7} (global convergence to the exact solution of the gradient
projection method). \emph{Suppose that assumptions of Theorem 5 hold and
also that parameters }$\lambda =\lambda \left( \delta ,3R\right) $\emph{\
and }$\beta =\beta \left( \delta ,3R\right) $\emph{\ are the same as in that
theorem. Let }$w_{n}=W_{n}+G,n=0,1,...$\emph{\ and }$w_{\min ,\lambda \left(
\delta ,3R\right) ,\beta \left( \delta ,3R\right) }=W_{\min ,\lambda \left(
\delta ,3R\right) ,\beta \left( \delta ,3R\right) }+G$ \emph{(Theorem 4).} 
\emph{Let }$c_{n,\text{comp}}\left( \mathbf{x}\right) $\emph{\ be the
function }$c\left( \mathbf{x}\right) $\emph{\ obtained from the function} $%
w_{n}\left( \mathbf{x},t\right) $ \emph{by the procedure outlined in the end
of section 3.} \emph{Then there exists a sufficiently small number }$\omega
_{1}=\omega _{1}\left( \Omega ,T,A,B,R,\gamma ,\lambda \right) \in \left(
0,\omega _{0}\right] $\emph{\ such that for any }$\omega \in \left( 0,\omega
_{1}\right) $\emph{\ there exists a number }$\theta =\theta \left( \omega
\right) \in \left( 0,1\right) $\emph{\ such that the the following
convergence estimates are valid for }$n=1,2,...$\emph{\ }%
\begin{equation}
\left\Vert w^{\ast }-w_{n}\right\Vert _{H^{2,1}\left( Q_{\gamma T}^{\pm
}\right) }\leq C_{2}\delta ^{\rho }+\theta ^{n}\left\Vert w_{\min ,\lambda
\left( \delta ,3R\right) ,\beta \left( \delta ,3R\right) }-w_{0}\right\Vert
_{H^{k_{n}}\left( Q_{T}^{\pm }\right) },\text{ }  \label{4.25}
\end{equation}%
\begin{equation}
\left\Vert c^{\ast }-c_{n,\text{comp}}\right\Vert _{L_{2}\left( \Omega
\right) }\leq C_{2}\delta ^{\rho }+\theta ^{n}\left\Vert w_{\min ,\lambda
\left( \delta ,3R\right) ,\beta \left( \delta ,3R\right) }-w_{0}\right\Vert
_{H^{k_{n}}\left( Q_{T}^{\pm }\right) }.  \label{4.26}
\end{equation}

\textbf{Remarks 3}:

\textbf{1.} \emph{Since the starting point }$W_{0}\in B_{0}\left( 2R\right) $%
\emph{\ of the iterative process (\ref{4.23}) is an arbitrary point of the
ball }$B_{0}\left( 2R\right) $\emph{\ and since }$R>0$\emph{\ is an
arbitrary number, then Theorem 7 ensures the \textbf{global convergence} of
the gradient projection method (\ref{4.23}) to the correct solution as long
as the noise level }$\delta $ \emph{tends to zero,} \emph{see section 1 for
our definition of the global convergence. }

\textbf{2.} \emph{We omit the proof of Theorem 3 below since, by Lemma 2.1
of \cite{Bak}, Theorem 3 follows immediately from Theorem 2. In addition, we
omit the proof of Theorem 6 since Theorem 2.1 of \cite{Bak} implies that
Theorem 6 follows immediately from Theorem 2.}

\section{Proof of Theorem 1}

\label{sec:5}

We prove this theorem only for functions $u\left( \mathbf{x},t\right) $ such
that 
\begin{equation}
u\in C^{3}\left( \overline{Q_{T}^{\pm }}\right) ,u\mid _{S_{T}^{\pm
}}=u_{x}\mid _{\Gamma _{T}^{\pm }}=0.  \label{5.70}
\end{equation}
The case $u\in H_{0}^{2,1}\left( Q_{T}^{\pm }\right) $ follows immediately
from this proof via density arguments. Below in this proof $O\left(
1/\lambda ^{k}\right) ,k\geq 1$ denotes different smooth functions, which
are independent on $u$ and for which the following estimate is valid $%
\left\Vert O\left( 1/\lambda ^{k}\right) \right\Vert _{C^{1}\left( \overline{%
Q_{T}^{\pm }}\right) }\leq C/\lambda ^{k},\forall \lambda \geq 1.$

Recall that by (\ref{3.10})\emph{\ }$\varphi _{\lambda }\left( x,t\right)
=\exp \left( 2\lambda \left( x^{2}-t^{2}\right) \right) .$ Introduce a new
function $v\left( \mathbf{x},t\right) =u\left( \mathbf{x},t\right) \exp
\left( \lambda \left( x^{2}-t^{2}\right) \right) .$ Then $u=v\exp \left(
-\lambda \left( x^{2}-t^{2}\right) \right) .$ Hence,%
\[
\left( u_{t}-\Delta u\right) =\left( v_{t}-\Delta v+4\lambda xv_{x}-4\lambda
^{2}x^{2}\left( 1-1/\left( 2\lambda x\right) \right) v+2\lambda tv\right)
\exp \left( -\lambda \left( x^{2}-t^{2}\right) \right) = 
\]%
\[
\left[ \left( -\Delta v-4\lambda ^{2}x^{2}\left( 1+O\left( 1/\lambda \right)
\right) v+2\lambda tv\right) +\left( v_{t}+4\lambda xv_{x}\right) \right]
\exp \left( -\lambda \left( x^{2}-t^{2}\right) \right) . 
\]%
Hence, 
\begin{equation}
\left( u_{t}-\Delta u\right) ^{2}\varphi _{\lambda }\geq \left(
2v_{t}+8\lambda xv_{x}\right) \left( -\Delta v-4\lambda ^{2}x^{2}\left(
1+O\left( 1/\lambda \right) \right) v+2\lambda tv\right) .  \label{5.1}
\end{equation}

\textbf{Step 1}. Estimate from the below the following term in (\ref{5.1}):

$2v_{t}\left( -\Delta v-4\lambda x^{2}\left( 1+O\left( 1/\lambda \right)
\right) v+2\lambda tv\right) ,$ 
\[
2v_{t}\left( -\Delta v-4\lambda x^{2}\left( 1+O\left( 1/\lambda \right)
\right) v+2\lambda tv\right) 
\]%
\[
=-2\mathop{\displaystyle \sum }\limits_{i=1}^{n}v_{x_{i}x_{i}}v_{t}+\left(
-4\lambda ^{2}x^{2}\left( 1+O\left( 1/\lambda \right) \right) v^{2}+2\lambda
tv^{2}\right) _{t}-2\lambda v^{2} 
\]%
\[
=\mathop{\displaystyle \sum }\limits_{i=1}^{n}\left( -2v_{x_{i}}v_{t}\right)
_{x_{i}}+2\mathop{\displaystyle \sum }\limits_{i=1}^{n}v_{x_{i}}v_{x_{i}t}+%
\left( -4\lambda ^{2}x^{2}\left( 1+O\left( 1/\lambda \right) \right)
v^{2}+2\lambda tv^{2}\right) _{t}-2\lambda v^{2} 
\]%
\[
=-2\lambda v^{2}+\mathop{\displaystyle \sum }\limits_{i=1}^{n}\left(
-2v_{x_{i}}v_{t}\right) _{x_{i}}+\left( \left( \nabla v\right) ^{2}-4\lambda
^{2}x^{2}\left( 1+O\left( 1/\lambda \right) \right) v^{2}+2\lambda
tv^{2}\right) _{t}. 
\]%
Thus, the desired estimate of Step 1 is:%
\begin{equation}
2v_{t}\left( -\Delta v-4\lambda x^{2}\left( 1+O\left( 1/\lambda \right)
\right) v+2\lambda tv\right) =-2\lambda v^{2}+\mathop{\rm div}\nolimits
U_{1}+V_{1t},  \label{5.2}
\end{equation}%
\begin{equation}
\mathop{\rm div}\nolimits U_{1}=\mathop{\displaystyle \sum }%
\limits_{i=1}^{n}\left( -2v_{x_{i}}v_{t}\right) _{x_{i}},  \label{5.3}
\end{equation}%
\begin{equation}
V_{1}=\left( \nabla v\right) ^{2}-4\lambda ^{2}x^{2}\left( 1+O\left(
1/\lambda \right) \right) v^{2}+2\lambda tv^{2}.  \label{5.4}
\end{equation}

\textbf{Step 2}. Estimate from the below the following term in (\ref{5.1}):

$8\lambda xv_{x}\left( -\Delta v-4\lambda ^{2}x^{2}\left( 1+O\left(
1/\lambda \right) \right) v+2\lambda tv\right) ,$

\[
8\lambda xv_{x}\left( -\Delta v-4\lambda ^{2}x^{2}\left( 1+O\left( 1/\lambda
\right) \right) v+2\lambda tv\right) 
\]%
\[
=-8\lambda xv_{x}v_{xx}+\mathop{\displaystyle \sum }\limits_{i=2}^{n}\left(
-8\lambda xv_{x}v_{x_{i}x_{i}}\right) 
\]%
\[
+\left( -16\lambda ^{3}x^{3}\left( 1+O\left( 1/\lambda \right) \right)
v^{2}+8\lambda ^{2}xtv^{2}\right) _{x}+48\lambda ^{3}x^{2}\left( 1+O\left(
1/\lambda \right) \right) v^{2} 
\]%
\[
=\left( -4\lambda xv_{x}^{2}\right) _{x}+4\lambda v_{x}^{2}+%
\mathop{\displaystyle \sum }\limits_{i=2}^{n}\left( -8\lambda
xv_{x}v_{x_{i}}\right) _{x_{i}}+\mathop{\displaystyle \sum }%
\limits_{i=2}^{n}\left( 8\lambda xv_{xx_{i}}v_{x_{i}}\right) 
\]%
\[
+48\lambda ^{3}x^{2}\left( 1+O\left( 1/\lambda \right) \right) v^{2}+\left(
-16\lambda ^{3}x^{3}\left( 1+O\left( 1/\lambda \right) \right)
v^{2}+8\lambda ^{2}xtv^{2}\right) _{x} 
\]%
\[
=4\lambda \left( v_{x}^{2}-\mathop{\displaystyle \sum }%
\limits_{i=2}^{n}v_{x_{i}}^{2}\right) +48\lambda ^{3}x^{2}\left( 1+O\left(
1/\lambda \right) \right) v^{2} 
\]%
\[
+\left( -4\lambda xv_{x}^{2}+4\lambda x\mathop{\displaystyle \sum }%
\limits_{i=2}^{n}v_{x_{i}}^{2}-16\lambda ^{3}x^{3}\left( 1+O\left( 1/\lambda
\right) \right) v^{2}+8\lambda ^{2}xtv^{2}\right) _{x} 
\]%
\[
+\mathop{\displaystyle \sum }\limits_{i=2}^{n}\left( -8\lambda
xv_{x}v_{x_{i}}\right) _{x_{i}}. 
\]

Thus, we end up with the following estimate of Step 2:%
\[
8\lambda xv_{x}\left( -\Delta v-4\lambda ^{2}x^{2}\left( 1+O\left( 1/\lambda
\right) \right) v+2\lambda tv\right) 
\]%
\begin{equation}
=4\lambda \left( v_{x}^{2}-\mathop{\displaystyle \sum }%
\limits_{i=2}^{n}v_{x_{i}}^{2}\right) +48\lambda ^{3}x^{2}\left( 1+O\left(
1/\lambda \right) \right) v^{2}+\mathop{\rm div}\nolimits U_{2},  \label{5.5}
\end{equation}%
\begin{equation}
\mathop{\rm div}\nolimits U_{2}=\left( -4\lambda xv_{x}^{2}+4\lambda x%
\mathop{\displaystyle \sum }\limits_{i=2}^{n}v_{x_{i}}^{2}-16\lambda
^{3}x^{3}\left( 1+O\left( 1/\lambda \right) \right) v^{2}+8\lambda
^{2}xtv^{2}\right) _{x}  \label{5.6}
\end{equation}%
\[
+\mathop{\displaystyle \sum }\limits_{i=2}^{n}\left( -8\lambda
xv_{x}v_{x_{i}}\right) _{x_{i}}. 
\]

\textbf{Step 3}. Analysis of boundary integrals over $S_{T}^{\pm }$.

Let $\nu =\nu \left( \mathbf{x}\right) $ be the unit outward looking normal
vector to $\partial \Omega $ at the point $x\in \partial \Omega .$ By Gauss'
formula, (\ref{5.3}) and (\ref{5.6}) 
\begin{equation}
\mathop{\displaystyle \int}\limits_{Q_{T}^{\pm }}\left( \mathop{\rm div}%
\nolimits U_{1}+\mathop{\rm div}\nolimits U_{2}\right) d\mathbf{x}dt=%
\mathop{\displaystyle \int}\limits_{-T}^{T}\mathop{\displaystyle \int}%
\limits_{\partial \Omega }\mathop{\displaystyle \sum }\limits_{i=1}^{n}%
\left( U_{1,i}+U_{2,i}\right) \cos \left( \nu \left( \mathbf{x}\right)
,x_{i}\right) dSdt,  \label{5.7}
\end{equation}%
where $U_{k}=\left( U_{k,1},...,U_{k,n}\right) ,k=1,2.$ Obviously, $%
U_{1i}=-2v_{x_{i}}v_{t}.$ Since (\ref{5.70}) holds and since 
\[
v_{t}\left( \mathbf{x},t\right) =\left( u_{t}-2\lambda tu\right) \left( 
\mathbf{x},t\right) \exp \left( \lambda \left( x^{2}-t^{2}\right) \right) , 
\]%
then $v_{t}\left( \mathbf{x},t\right) =0$ for $\mathbf{x}\in \partial \Omega
.$ Hence, in (\ref{5.7})%
\begin{equation}
\mathop{\displaystyle \int}\limits_{-T}^{T}\mathop{\displaystyle \int}%
\limits_{\partial \Omega }\mathop{\displaystyle \sum }%
\limits_{i=1}^{n}U_{1,i}\cos \left( \nu \left( \mathbf{x}\right)
,x_{i}\right) dSdt=0.  \label{5.8}
\end{equation}%
We now analyze the first term in the right hand side of (\ref{5.6}). We have%
\begin{equation}
v_{x}\left( \mathbf{x},t\right) =\left( u_{x}+2\lambda xu\right) \left( 
\mathbf{x},t\right) \exp \left( \lambda \left( x^{2}-t^{2}\right) \right) ,
\label{5.9}
\end{equation}%
\begin{equation}
v_{x_{i}}\left( \mathbf{x},t\right) =u_{x_{i}}\left( \mathbf{x},t\right)
\exp \left( \lambda \left( x^{2}-t^{2}\right) \right) ,\text{ }i=2,...,n.
\label{5.10}
\end{equation}%
By (\ref{5.6}), (\ref{5.70}), (\ref{5.9}) and (\ref{5.10}) 
\begin{equation}
\mathop{\displaystyle \int}\limits_{-T}^{T}\mathop{\displaystyle \int}%
\limits_{\partial \Omega }U_{2,1}\cos \left( \nu \left( \mathbf{x}\right)
,x_{1}\right) dSdt=4\lambda A\mathop{\displaystyle \int}\limits_{-T}^{T}%
\mathop{\displaystyle \int}\limits_{\Gamma ^{\prime }}u_{x}^{2}\left( A,%
\overline{x}\right) dS\geq 0,  \label{5.11}
\end{equation}%
where $\Gamma ^{\prime }=\left\{ x=A\right\} \cap \partial \Omega .$
Similarly%
\begin{equation}
\mathop{\displaystyle \int}\limits_{-T}^{T}\mathop{\displaystyle \int}%
\limits_{\partial \Omega }\mathop{\displaystyle \sum }%
\limits_{i=2}^{n}U_{2,i}\cos \left( \nu \left( \mathbf{x}\right)
,x_{i}\right) dSdt=0.  \label{5.12}
\end{equation}%
Using (\ref{5.7}), (\ref{5.8}), (\ref{5.11}) and (\ref{5.12}), we obtain%
\begin{equation}
\mathop{\displaystyle \int}\limits_{Q_{T}^{\pm }}\left( \mathop{\rm div}%
\nolimits U_{1}+\mathop{\rm div}\nolimits U_{2}\right) d\mathbf{x}dt\geq 0.
\label{5.13}
\end{equation}

\textbf{Step 4}. Integrate (\ref{5.1}) over $Q_{T}^{\pm }.$ Then sum up (\ref%
{5.2}) with (\ref{5.5}), integrate the resulting inequality over $Q_{T}^{\pm
}$ and use that integral of (\ref{5.1}), Gauss' formula, (\ref{5.4}) and (%
\ref{5.8})-(\ref{5.12}). We obtain for all $\lambda \geq \lambda _{0}$ and
all $u\in C^{2}\left( \overline{Q_{T}^{\pm }}\right) \cap H_{0}^{2,1}\left(
Q_{T}^{\pm }\right) $%
\begin{equation}
\mathop{\displaystyle \int}\limits_{Q_{T}^{\pm }}\left( u_{t}-\Delta
u\right) ^{2}\varphi _{\lambda }d\mathbf{x}dt\geq -4\lambda %
\mathop{\displaystyle \int}\limits_{Q_{T}^{\pm }}\left( \nabla u\right)
^{2}\varphi _{\lambda }d\mathbf{x}dt+47\lambda ^{3}\mathop{\displaystyle
\int}\limits_{Q_{T}^{\pm }}u^{2}x^{2}\varphi _{\lambda }d\mathbf{x}dt
\label{5.14}
\end{equation}%
\[
-C\exp \left( 2\lambda \left( B^{2}-T^{2}\right) \right) \mathop{%
\displaystyle \int}\limits_{\Omega }\left[ \left( \left( \nabla u\right)
^{2}+\lambda ^{2}u^{2}\right) \left( \mathbf{x},T\right) +\left( \left(
\nabla u\right) ^{2}+\lambda ^{2}u^{2}\right) \left( \mathbf{x},-T\right) %
\right] d\mathbf{x.} 
\]%
The inconvenient point of (\ref{5.14}) is the presence of the negative term
in the first line of (\ref{5.14}). Therefore, we continue.

\textbf{Step 5}. Estimate from the below $\left( u_{t}-\Delta u\right)
u\varphi _{\lambda },$ and then estimate the corresponding integral over $%
Q_{T}^{\pm }.$%
\[
\left( u_{t}-\Delta u\right) u\varphi _{\lambda }=\left( \frac{u^{2}}{2}%
\varphi _{\lambda }\right) _{t}+2\lambda u^{2}\varphi _{\lambda }+\left(
-u_{x}u\varphi _{\lambda }\right) _{x}+u_{x}^{2}\varphi _{\lambda }+4\lambda
xu_{x}u\varphi _{\lambda } 
\]%
\[
+\mathop{\displaystyle \sum }\limits_{i=2}^{n}\left( -u_{x_{i}}u\varphi
_{\lambda }\right) _{x_{i}}+\mathop{\displaystyle \sum }%
\limits_{i=2}^{n}u_{x_{i}}^{2}\varphi _{\lambda } 
\]%
\[
=\left( \nabla u\right) ^{2}\varphi _{\lambda }+\mathop{\displaystyle \sum }%
\limits_{i=1}^{n}\left( -u_{x_{i}}u\varphi _{\lambda }\right)
_{x_{i}}+\left( \frac{u^{2}}{2}\varphi _{\lambda }\right) _{t}+\left(
2\lambda xu^{2}\varphi _{\lambda }\right) _{x} 
\]%
\[
-8\lambda ^{2}x^{2}\left( 1+O\left( 1/\lambda \right) \right) u^{2}\varphi
_{\lambda }. 
\]%
Hence,%
\begin{equation}
\left( u_{t}-\Delta u\right) u\varphi _{\lambda }\geq \left( \nabla u\right)
^{2}\varphi _{\lambda }-9\lambda ^{2}x^{2}u^{2}\varphi _{\lambda }+%
\mathop{\rm div}\nolimits U_{3}+V_{2t},  \label{5.15}
\end{equation}%
\begin{equation}
\mathop{\rm div}\nolimits U_{3}=\mathop{\displaystyle \sum }%
\limits_{i=1}^{n}\left( -u_{x_{i}}u\varphi _{\lambda }\right)
_{x_{i}}+\left( 2\lambda xu^{2}\varphi _{\lambda }\right) _{x},  \label{5.16}
\end{equation}%
\begin{equation}
V_{2}=\frac{u^{2}}{2}\varphi _{\lambda }.  \label{5.17}
\end{equation}%
Hence, by (\ref{5.70}), (\ref{5.16}) and Gauss formula 
\begin{equation}
\mathop{\displaystyle \int}\limits_{Q_{T}^{\pm }}\mathop{\rm div}\nolimits
U_{3}d\mathbf{x}dt=0.  \label{5.18}
\end{equation}%
Integrate (\ref{5.15}) over $Q_{T}^{\pm }$ using (\ref{5.17}) and (\ref{5.18}%
). Then multiply the resulting inequality by $5\lambda $ and sum up with (%
\ref{5.14}). We obtain%
\[
5\lambda \mathop{\displaystyle \int}\limits_{Q_{T}^{\pm }}\left(
u_{t}-\Delta u\right) u\varphi _{\lambda }dxdt+\mathop{\displaystyle \int}%
\limits_{Q_{T}^{\pm }}\left( u_{t}-\Delta u\right) ^{2}\varphi _{\lambda }d%
\mathbf{x}dt 
\]%
\begin{equation}
\geq \lambda \mathop{\displaystyle \int}\limits_{Q_{T}^{\pm }}\left( \nabla
u\right) ^{2}\varphi _{\lambda }d\mathbf{x}dt+2\lambda ^{3}%
\mathop{\displaystyle \int}\limits_{Q_{T}^{\pm }}u^{2}x^{2}\varphi _{\lambda
}d\mathbf{x}dt  \label{5.19}
\end{equation}%
\[
-C\exp \left( 2\lambda \left( B^{2}-T^{2}\right) \right) \mathop{%
\displaystyle \int}\limits_{\Omega }\left[ \left( \left( \nabla u\right)
^{2}+\lambda ^{2}u^{2}\right) \left( \mathbf{x},T\right) +\left( \left(
\nabla u\right) ^{2}+\lambda ^{2}u^{2}\right) \left( \mathbf{x},-T\right) %
\right] d\mathbf{x.} 
\]%
Next, by the Cauchy-Schwarz inequality%
\begin{equation}
5\lambda \mathop{\displaystyle \int}\limits_{Q_{T}^{\pm }}\left(
u_{t}-\Delta u\right) u\varphi _{\lambda }dxdt+\mathop{\displaystyle \int}%
\limits_{Q_{T}^{\pm }}\left( u_{t}-\Delta u\right) ^{2}\varphi _{\lambda }d%
\mathbf{x}dt  \label{5.20}
\end{equation}%
\[
\leq \frac{7}{2}\mathop{\displaystyle \int}\limits_{Q_{T}^{\pm }}\left(
u_{t}-\Delta u\right) ^{2}\varphi _{\lambda }d\mathbf{x}dt+\frac{5}{2}%
\lambda ^{2}\mathop{\displaystyle \int}\limits_{Q_{T}^{\pm }}u^{2}\varphi
_{\lambda }d\mathbf{x}dt. 
\]%
Since for sufficiently large $\lambda _{0}>1$ and for $\lambda \geq \lambda
_{0}$ 
\begin{equation}
2\lambda ^{3}\mathop{\displaystyle \int}\limits_{Q_{T}^{\pm
}}u^{2}x^{2}\varphi _{\lambda }d\mathbf{x}dt-\frac{5}{2}\lambda ^{2}%
\mathop{\displaystyle \int}\limits_{Q_{T}^{\pm }}u^{2}\varphi _{\lambda }d%
\mathbf{x}dt\geq \lambda ^{3}\mathop{\displaystyle \int}\limits_{Q_{T}^{\pm
}}u^{2}x^{2}\varphi _{\lambda }d\mathbf{x}dt,  \label{5.21}
\end{equation}%
then (\ref{5.19})-(\ref{5.21}) imply that for all $u\in H_{0}^{2,1}\left(
Q_{T}^{\pm }\right) $ and for all $\lambda \geq \lambda _{0}$ 
\begin{equation}
\mathop{\displaystyle \int}\limits_{Q_{T}^{\pm }}\left( u_{t}-\Delta
u\right) ^{2}\varphi _{\lambda }d\mathbf{x}dt\geq C\lambda %
\mathop{\displaystyle \int}\limits_{Q_{T}^{\pm }}\left[ \left( \nabla
u\right) ^{2}+\lambda ^{2}u^{2}\right] \varphi _{\lambda }d\mathbf{x}dt
\label{5.22}
\end{equation}%
\[
-C\exp \left( 2\lambda \left( B^{2}-T^{2}\right) \right) \mathop{%
\displaystyle \int}\limits_{\Omega }\left[ \left( \left( \nabla u\right)
^{2}+\lambda ^{2}u^{2}\right) \left( \mathbf{x},T\right) +\left( \left(
\nabla u\right) ^{2}+\lambda ^{2}u^{2}\right) \left( \mathbf{x},-T\right) %
\right] d\mathbf{x,} 
\]

which is a part of estimate (\ref{4.2}). We now need to incorporate in our
estimate terms with $u_{t}^{2},u_{x_{i}x_{j}}^{2}.$

\textbf{Step 6}. Incorporating terms with $u_{t}^{2},u_{x_{i}x_{j}}^{2}.$

We have%
\begin{equation}
\left( u_{t}-\Delta u\right) ^{2}\varphi _{\lambda }=u_{t}^{2}\varphi
_{\lambda }-2u_{t}u_{xx}\varphi _{\lambda }-\mathop{\displaystyle \sum }%
\limits_{i=2}^{n}2u_{t}u_{x_{i}x_{i}}\varphi _{\lambda }+\left( \Delta
u\right) ^{2}\varphi _{\lambda }.  \label{5.23}
\end{equation}%
Denote 
\begin{equation}
z_{1}=-2u_{t}u_{xx}\varphi _{\lambda },\text{ }z_{2}=-\mathop{\displaystyle
\sum }\limits_{i=2}^{n}2u_{t}u_{x_{i}x_{i}}\varphi _{\lambda },\text{ }%
z_{3}=\left( \Delta u\right) ^{2}\varphi _{\lambda }  \label{5.24}
\end{equation}%
and estimate each of terms in (\ref{5.24}). First, we have%
\[
z_{1}=-2u_{t}u_{xx}\varphi _{\lambda }=\left( -2u_{t}u_{x}\varphi _{\lambda
}\right) _{x}+2u_{tx}u_{x}\varphi _{\lambda }+8\lambda xu_{t}u_{x}\varphi
_{\lambda } 
\]%
\[
=\left( -2u_{t}u_{x}\varphi _{\lambda }\right) _{x}+\left( u_{x}^{2}\varphi
_{\lambda }\right) _{t}+2\lambda tu_{x}^{2}\varphi _{\lambda }+8\lambda
xu_{t}u_{x}\varphi _{\lambda } 
\]%
\[
\geq -\frac{1}{2}u_{t}^{2}\varphi _{\lambda }-C\lambda ^{2}u_{x}^{2}\varphi
_{\lambda }+\left( -2u_{t}u_{x}\varphi _{\lambda }\right) _{x}+\left(
u_{x}^{2}\varphi _{\lambda }\right) _{t}, 
\]%
Thus,%
\begin{equation}
z_{1}\geq -\frac{1}{2}u_{t}^{2}\varphi _{\lambda }-C\lambda
^{2}u_{x}^{2}\varphi _{\lambda }+\left( -2u_{t}u_{x}\varphi _{\lambda
}\right) _{x}+\left( u_{x}^{2}\varphi _{\lambda }\right) _{t}.  \label{5.25}
\end{equation}%
We now estimate $z_{2},$%
\[
z_{2}=\mathop{\displaystyle \sum }\limits_{i=2}^{n}\left(
-2u_{t}u_{x_{i}}\varphi _{\lambda }\right) _{x_{i}}+%
\mathop{\displaystyle
\sum }\limits_{i=2}^{n}2u_{tx_{i}}u_{x_{i}}\varphi _{\lambda } 
\]%
\[
=\mathop{\displaystyle \sum }\limits_{i=2}^{n}\left( u_{x_{i}}^{2}\varphi
_{\lambda }\right) _{t}+4\lambda t\mathop{\displaystyle \sum }%
\limits_{i=2}^{n}u_{x_{i}}^{2}\varphi _{\lambda }+%
\mathop{\displaystyle \sum
}\limits_{i=2}^{n}\left( -2u_{t}u_{x_{i}}\varphi _{\lambda }\right) _{x_{i}} 
\]%
\[
\geq -C\lambda \mathop{\displaystyle \sum }\limits_{i=2}^{n}u_{x_{i}}^{2}%
\varphi _{\lambda }+\mathop{\displaystyle \sum }\limits_{i=2}^{n}\left(
u_{x_{i}}^{2}\varphi _{\lambda }\right) _{t}+\mathop{\displaystyle \sum }%
\limits_{i=2}^{n}\left( -2u_{t}u_{x_{i}}\varphi _{\lambda }\right) _{x_{i}}. 
\]%
Thus, 
\begin{equation}
z_{2}\geq -C\lambda \mathop{\displaystyle \sum }%
\limits_{i=2}^{n}u_{x_{i}}^{2}\varphi _{\lambda }+\mathop{\displaystyle \sum
}\limits_{i=2}^{n}\left( u_{x_{i}}^{2}\varphi _{\lambda }\right) _{t}+%
\mathop{\displaystyle \sum }\limits_{i=2}^{n}\left( -2u_{t}u_{x_{i}}\varphi
_{\lambda }\right) _{x_{i}}.  \label{5.26}
\end{equation}

Now we estimate $z_{3},$%
\[
z_{3}=\left( \Delta u\right) ^{2}\varphi _{\lambda }=\left( u_{xx}+%
\mathop{\displaystyle \sum }\limits_{i=2}^{n}u_{x_{i}x_{i}}\right)
^{2}\varphi _{\lambda } 
\]%
\[
=\mathop{\displaystyle \sum }\limits_{i=1}^{n}u_{x_{i}x_{i}}^{2}\varphi
_{\lambda }+2\mathop{\displaystyle \sum }%
\limits_{i=2}^{n}u_{xx}u_{x_{i}x_{i}}\varphi _{\lambda }+2%
\mathop{\displaystyle \sum }\limits_{i,j=2,i\neq
j}^{n}u_{x_{i}x_{i}}u_{x_{j}x_{j}}\varphi _{\lambda } 
\]%
\[
=\mathop{\displaystyle \sum }\limits_{i=1}^{n}u_{x_{i}x_{i}}^{2}\varphi
_{\lambda }+\left( 2\mathop{\displaystyle \sum }%
\limits_{i=2}^{n}u_{x}u_{x_{i}x_{i}}\varphi _{\lambda }\right) _{x}-8\lambda
x\mathop{\displaystyle \sum }\limits_{i=2}^{n}u_{x}u_{x_{i}x_{i}}\varphi
_{\lambda } 
\]%
\begin{equation}
-2\mathop{\displaystyle \sum }\limits_{i=2}^{n}u_{x}u_{xx_{i}x_{i}}\varphi
_{\lambda }+\left( 2\mathop{\displaystyle \sum }\limits_{i,j=2,i\neq
j}^{n}u_{x_{i}}u_{x_{j}x_{j}}\varphi _{\lambda }\right) _{x_{i}}-2%
\mathop{\displaystyle \sum }\limits_{i,j=2,i\neq
j}^{n}u_{x_{i}}u_{x_{i}x_{j}x_{j}}\varphi _{\lambda }  \label{5.27}
\end{equation}%
\[
=\mathop{\displaystyle \sum }\limits_{i=1}^{n}u_{x_{i}x_{i}}^{2}\varphi
_{\lambda }+\left( 2\mathop{\displaystyle \sum }%
\limits_{i=2}^{n}u_{x}u_{x_{i}x_{i}}\varphi _{\lambda }\right) _{x}-8\lambda
x\mathop{\displaystyle \sum }\limits_{i=2}^{n}u_{x}u_{x_{i}x_{i}}\varphi
_{\lambda } 
\]%
\[
+\left( -2\mathop{\displaystyle \sum }\limits_{i=2}^{n}u_{x}u_{xx_{i}}%
\varphi _{\lambda }\right) _{x_{i}}+\mathop{\displaystyle \sum }%
\limits_{i=2}^{n}u_{xx_{i}}^{2}\varphi _{\lambda }+\left( 2%
\mathop{\displaystyle \sum }\limits_{i,j=2,i\neq
j}^{n}u_{x_{i}}u_{x_{j}x_{j}}\varphi _{\lambda }\right) _{x_{i}} 
\]%
\[
+\left( -2\mathop{\displaystyle \sum }\limits_{i,j=2,i\neq
j}^{n}u_{x_{i}}u_{x_{i}x_{j}}\varphi _{\lambda }\right) _{x_{j}}+2%
\mathop{\displaystyle \sum }\limits_{i,j=2,i\neq
j}^{n}u_{x_{i}x_{j}}^{2}\varphi _{\lambda }. 
\]%
Since by the Cauchy-Schwarz inequality%
\[
-8\lambda x\mathop{\displaystyle \sum }\limits_{i=2}^{n}u_{x}u_{x_{i}x_{i}}%
\varphi _{\lambda }\geq -C\lambda ^{2}\left( \nabla u\right) ^{2}\varphi
_{\lambda }-\frac{1}{2}\mathop{\displaystyle \sum }%
\limits_{i=2}^{n}u_{x_{i}x_{i}}^{2}\varphi _{\lambda }, 
\]%
then (\ref{5.27}) implies that%
\[
z_{3}\geq \frac{1}{2}\mathop{\displaystyle \sum }%
\limits_{i,j=1}^{n}u_{x_{i}x_{j}}^{2}\varphi _{\lambda }-C\lambda ^{2}\left(
\nabla u\right) ^{2}\varphi _{\lambda } 
\]%
\begin{equation}
+\left( 2\mathop{\displaystyle \sum }\limits_{i=2}^{n}u_{x}u_{x_{i}x_{i}}%
\varphi _{\lambda }\right) _{x}+\left( -2\mathop{\displaystyle \sum }%
\limits_{i=2}^{n}u_{x}u_{xx_{i}}\varphi _{\lambda }\right) _{x_{i}}
\label{5.28}
\end{equation}%
\[
+\left( -2\mathop{\displaystyle \sum }\limits_{i,j=2,i\neq
j}^{n}u_{x_{i}}u_{x_{i}x_{j}}\varphi _{\lambda }\right) _{x_{j}}+\left( 2%
\mathop{\displaystyle \sum }\limits_{i,j=2,i\neq
j}^{n}u_{x_{i}}u_{x_{j}x_{j}}\varphi _{\lambda }\right) _{x_{i}}. 
\]%
Combining (\ref{5.23})-(\ref{5.28}), we obtain%
\begin{equation}
\frac{1}{4\lambda }\left( u_{t}-\Delta u\right) ^{2}\varphi _{\lambda }\geq 
\frac{1}{8\lambda }\left( u_{t}^{2}+\mathop{\displaystyle \sum }%
\limits_{i,j=1}^{n}u_{x_{i}x_{j}}^{2}\right) \varphi _{\lambda }-\frac{C}{2}%
\lambda \left( \nabla u\right) ^{2}\varphi _{\lambda }  \label{5.29}
\end{equation}%
\[
+\mathop{\rm div}\nolimits U_{4}+\left( u_{x}^{2}\varphi _{\lambda }\right)
_{t}, 
\]%
\begin{equation}
\mathop{\displaystyle \int}\limits_{Q_{T}^{\pm }}\mathop{\rm div}\nolimits
U_{4}d\mathbf{x}dt=0.  \label{5.30}
\end{equation}

Using (\ref{5.30}), integrate (\ref{5.29}) over $Q_{T}^{\pm }.$ Then sum up
the resulting inequality with (\ref{5.22}). Then we obtain the target
estimate (\ref{4.2}) of this theorem. $\ \square $

\section{Proofs of Theorems 2 and 4}

\label{sec:6}

Lemma 1 follows immediately either from Lemma 1.10.3 of \cite{BKbook}\emph{\ 
}or from Lemma 3.1 of \cite{Ksurvey}.

\textbf{Lemma 1}. \emph{The following estimate holds for every function }$%
q\in L_{2}\left( Q_{T}^{\pm }\right) $\emph{\ and for every }$\lambda \geq
1: $%
\[
\mathop{\displaystyle \int}\limits_{Q_{T}^{\pm }}\left( \mathop{%
\displaystyle \int}\limits_{0}^{t}q\left( \mathbf{x},\tau \right) d\tau
\right) ^{2}\varphi _{\lambda }\left( x,t\right) d\mathbf{x}dt\leq \frac{1}{%
4\lambda }\mathop{\displaystyle \int}\limits_{Q_{T}^{\pm }}q^{2}\left( 
\mathbf{x},t\right) \varphi _{\lambda }\left( x,t\right) d\mathbf{x}dt. 
\]

\subsection{Proof of Theorem 2}

\label{sec:6.1}

Let $w_{1},w_{2}\in \overline{B\left( R,p_{0},p_{1}\right) }$ be two
arbitrary functions. Denote $h=w_{2}-w_{1}.$ Then $w_{2}=w_{1}+h$ and also 
\begin{equation}
h\in \overline{B_{0}\left( 2R\right) }.  \label{6.1}
\end{equation}%
First, we evaluate the expression $\left( K\left( w_{1}+h\right) \right)
^{2}-\left( K\left( w_{1}\right) \right) ^{2},$ where the nonlinear operator 
$K$ is given in (\ref{3.8}). We have%
\[
\left( K\left( w_{1}+h\right) \right) ^{2}= 
\]%
\[
\left[ h_{t}-Lh-2\nabla h\mathop{\displaystyle \int}\limits_{0}^{t}\nabla
w_{1}\left( \mathbf{x},\tau \right) d\tau -2\nabla w_{1}\mathop{%
\displaystyle \int}\limits_{0}^{t}\nabla h_{1}\left( \mathbf{x},\tau \right)
d\tau -2\nabla h\mathop{\displaystyle \int}\limits_{0}^{t}\nabla h\left( 
\mathbf{x},\tau \right) d\tau +K\left( w_{1}\right) \right] ^{2} 
\]%
\[
=\left( h_{t}-Lh-2\nabla h\mathop{\displaystyle \int}\limits_{0}^{t}\nabla
w_{1}\left( \mathbf{x},\tau \right) d\tau -2\nabla w_{1}\mathop{%
\displaystyle \int}\limits_{0}^{t}\nabla h_{1}\left( \mathbf{x},\tau \right)
d\tau -2\nabla h\mathop{\displaystyle \int}\limits_{0}^{t}\nabla h\left( 
\mathbf{x},\tau \right) d\tau \right) ^{2} 
\]%
\[
-4K\left( w_{1}\right) \nabla h\mathop{\displaystyle \int}%
\limits_{0}^{t}\nabla h\left( \mathbf{x},\tau \right) d\tau +\left( K\left(
w_{1}\right) \right) ^{2} 
\]%
\[
+2K\left( w_{1}\right) \left( h_{t}-Lh-2\nabla h\mathop{\displaystyle \int}%
\limits_{0}^{t}\nabla w_{1}\left( \mathbf{x},\tau \right) d\tau -2\nabla
w_{1}\mathop{\displaystyle \int}\limits_{0}^{t}\nabla h_{1}\left( \mathbf{x}%
,\tau \right) d\tau \right) . 
\]%
Let $Lin\left( h\right) $ be the linear, with respect to $h$, part of the
above expression,%
\begin{equation}
Lin\left( h\right) =2K\left( w_{1}\right) \left( h_{t}-Lh-2\nabla h%
\mathop{\displaystyle \int}\limits_{0}^{t}\nabla w_{1}\left( \mathbf{x},\tau
\right) d\tau -2\nabla w_{1}\mathop{\displaystyle \int}\limits_{0}^{t}\nabla
h_{1}\left( \mathbf{x},\tau \right) d\tau \right) .  \label{6.2}
\end{equation}%
Then%
\begin{equation}
\left( K\left( w_{1}+h\right) \right) ^{2}-\left( K\left( w_{1}\right)
\right) ^{2}=Lin\left( h\right)  \label{6.3}
\end{equation}%
\[
+\left( h_{t}-Lh-2\nabla h\mathop{\displaystyle \int}\limits_{0}^{t}\nabla
w_{1}\left( \mathbf{x},\tau \right) d\tau -2\nabla w_{1}\mathop{%
\displaystyle \int}\limits_{0}^{t}\nabla h_{1}\left( \mathbf{x},\tau \right)
d\tau -2\nabla h\mathop{\displaystyle \int}\limits_{0}^{t}\nabla h\left( 
\mathbf{x},\tau \right) d\tau \right) ^{2} 
\]%
\[
-4K\left( w_{1}\right) \nabla h\mathop{\displaystyle \int}%
\limits_{0}^{t}\nabla h\left( \mathbf{x},\tau \right) d\tau . 
\]%
Using (\ref{3.12}), (\ref{3.13}), (\ref{6.1}) and the Cauchy-Schwarz
inequality, we obtain 
\[
\left( h_{t}-Lh-2\nabla h\mathop{\displaystyle \int}\limits_{0}^{t}\nabla
w_{1}\left( \mathbf{x},\tau \right) d\tau -2\nabla w_{1}\mathop{%
\displaystyle \int}\limits_{0}^{t}\nabla h_{1}\left( \mathbf{x},\tau \right)
d\tau -2\nabla h\mathop{\displaystyle \int}\limits_{0}^{t}\nabla h\left( 
\mathbf{x},\tau \right) d\tau \right) ^{2} 
\]%
\begin{equation}
-4K\left( w_{1}\right) \nabla h\mathop{\displaystyle \int}%
\limits_{0}^{t}\nabla h\left( \mathbf{x},\tau \right) d\tau  \label{6.4}
\end{equation}%
\[
\geq \frac{1}{2}\left( h_{t}-Lh\right) ^{2}-C_{1}\left( \nabla h\right)
^{2}-C_{1}\left( \mathop{\displaystyle \int}\limits_{0}^{t}\nabla h\left( 
\mathbf{x},\tau \right) d\tau \right) ^{2}. 
\]%
Hence, (\ref{3.15}) and (\ref{6.2})-(\ref{6.4}) lead to%
\[
J_{\lambda ,\beta }\left( w_{1}+h\right) -J_{\lambda ,\beta }\left(
w_{1}\right) -e^{-2\lambda B^{2}}\mathop{\displaystyle \int}%
\limits_{Q_{T}^{\pm }}Lin\left( h\right) \varphi _{\lambda }d\mathbf{x}%
dt+2\beta \left\{ w,h\right\} 
\]%
\begin{equation}
\geq \frac{1}{2}e^{-2\lambda B^{2}}\mathop{\displaystyle \int}%
\limits_{Q_{T}^{\pm }}\left( h_{t}-Lh\right) ^{2}\varphi _{\lambda }d\mathbf{%
x}dt  \label{6.5}
\end{equation}%
\[
-C_{1}e^{-2\lambda B^{2}}\mathop{\displaystyle \int}\limits_{Q_{T}^{\pm }}%
\left[ \left( \nabla h\right) ^{2}+\left( \mathop{\displaystyle \int}%
\limits_{0}^{t}\nabla h\left( \mathbf{x},\tau \right) d\tau \right) ^{2}%
\right] \varphi _{\lambda }d\mathbf{x}dt+\beta \left\Vert h\right\Vert
_{H^{k_{n}}\left( Q_{T}^{\pm }\right) }^{2}. 
\]%
Here and below $\left\{ ,\right\} $ is the scalar product in $%
H^{k_{n}}\left( Q_{T}^{\pm }\right) .$

Consider now the functional $S\left( h\right) :H_{0}^{k_{n}}\left(
Q_{T}^{\pm }\right) \rightarrow \mathbb{R}$ defined as%
\begin{equation}
S\left( h\right) =e^{-2\lambda B^{2}}\mathop{\displaystyle \int}%
\limits_{Q_{T}^{\pm }}Lin\left( h\right) \varphi _{\lambda }d\mathbf{x}%
dt+2\beta \left\{ w,h\right\} .  \label{6.6}
\end{equation}%
It is clear from (\ref{6.2}) that $S\left( h\right) $ is a bounded linear
functional. Hence, by Riesz theorem there exists a function $Z\in
H_{0}^{k_{n}}\left( Q_{T}^{\pm }\right) $ such that $S\left( h\right)
=\left\{ Z,h\right\} ,\forall h\in H_{0}^{k_{n}}\left( Q_{T}^{\pm }\right) .$
Furthermore, it follows from (\ref{6.3}) that 
\[
\lim_{\left\Vert h\right\Vert _{H^{k_{n}}\left( Q_{T}^{\pm }\right)
}\rightarrow 0}\left\vert J_{\lambda ,\beta }\left( w_{1}+h\right)
-J_{\lambda ,\beta }\left( w_{1}\right) -S\left( h\right) \right\vert =0. 
\]%
Hence, $S\left( h\right) $ is the Frech\'{e}t derivative of the functional $%
J_{\lambda ,\beta }\left( w\right) $ at the point $w_{1},$%
\begin{equation}
S\left( h\right) =J_{\lambda ,\beta }^{\prime }\left( w_{1}\right) \left(
h\right) =\left\{ Z,h\right\} =\left\{ J_{\lambda ,\beta }^{\prime }\left(
w_{1}\right) ,h\right\} ,\forall h\in H_{0}^{k_{n}}\left( Q_{T}^{\pm
}\right) ,  \label{6.7}
\end{equation}%
i.e. we can set $Z=J_{\lambda ,\beta }^{\prime }\left( w_{1}\right) .$ Note
that the proof of the existence of the Frech\'{e}t derivative on the set $%
B\left( 3R,p_{0},p_{1}\right) $, as claimed in this theorem, is basically
the same as the one above. Thus, (\ref{6.5})-(\ref{6.7}) imply that 
\[
J_{\lambda ,\beta }\left( w_{1}+h\right) -J_{\lambda ,\beta }\left(
w_{1}\right) -J_{\lambda ,\beta }^{\prime }\left( w_{1}\right) \left(
h\right) 
\]%
\begin{equation}
\geq \frac{1}{2}e^{-2\lambda B^{2}}\mathop{\displaystyle \int}%
\limits_{Q_{T}^{\pm }}\left( h_{t}-Lh\right) ^{2}\varphi _{\lambda }d\mathbf{%
x}dt  \label{6.8}
\end{equation}%
\[
-C_{1}e^{-2\lambda B^{2}}\mathop{\displaystyle \int}\limits_{Q_{T}^{\pm }}%
\left[ \left( \nabla h\right) ^{2}+\left( \mathop{\displaystyle \int}%
\limits_{0}^{t}\nabla h\left( \mathbf{x},\tau \right) d\tau \right) ^{2}%
\right] \varphi _{\lambda }d\mathbf{x}dt+\beta \left\Vert h\right\Vert
_{H^{k_{n}}\left( Q_{T}^{\pm }\right) }^{2}. 
\]%
Applying Theorem 1 and Lemma 1, we obtain for all $\lambda \geq \lambda
_{0}\geq 1$ 
\[
\frac{1}{2}e^{-2\lambda B^{2}}\mathop{\displaystyle \int}\limits_{Q_{T}^{\pm
}}\left( h_{t}-Lh\right) ^{2}\varphi _{\lambda }d\mathbf{x}dt 
\]%
\begin{equation}
-C_{1}e^{-2\lambda B^{2}}\mathop{\displaystyle \int}\limits_{Q_{T}^{\pm }}%
\left[ \left( \nabla h\right) ^{2}+\left( \mathop{\displaystyle \int}%
\limits_{0}^{t}\nabla h\left( \mathbf{x},\tau \right) d\tau \right) ^{2}%
\right] \varphi _{\lambda }d\mathbf{x}dt+\beta \left\Vert h\right\Vert
_{H^{k_{n}}\left( Q_{T}^{\pm }\right) }^{2}  \label{6.9}
\end{equation}%
\[
\geq \frac{C}{\lambda }e^{-2\lambda B^{2}}\mathop{\displaystyle \int}%
\limits_{Q_{T}^{\pm }}\left( h_{t}^{2}+\mathop{\displaystyle \sum }%
\limits_{i,j=1}^{n}h_{x_{i}x_{j}}^{2}\right) \varphi _{\lambda }d\mathbf{x}%
dt 
\]%
\[
+C\lambda e^{-2\lambda B^{2}}\mathop{\displaystyle \int}\limits_{Q_{T}^{\pm
}}\left[ \left( \nabla h\right) ^{2}+\lambda ^{2}h^{2}\right] \varphi
_{\lambda }d\mathbf{x}dt-C_{1}e^{-2\lambda B^{2}}\mathop{\displaystyle \int}%
\limits_{Q_{T}^{\pm }}\left( \nabla h\right) ^{2}\varphi _{\lambda }d\mathbf{%
x}dt 
\]%
\[
-Ce^{-2\lambda T^{2}}\mathop{\displaystyle \int}\limits_{\Omega }\left[
\left( h_{t}^{2}+\left( \nabla h\right) ^{2}+\lambda ^{2}h^{2}\right) \left( 
\mathbf{x},T\right) +\left( h_{t}^{2}+\left( \nabla h\right) ^{2}+\lambda
^{2}h^{2}\right) \left( \mathbf{x},-T\right) \right] d\mathbf{x}\text{ } 
\]%
\[
+\beta \left\Vert h\right\Vert _{H^{k_{n}}\left( Q_{T}^{\pm }\right) }^{2}. 
\]%
Choose $\lambda _{1}\geq \lambda _{0}$ so large that $C\lambda \geq 2C_{1}$
and also $2e^{-\lambda T^{2}}\geq C\lambda ^{2}e^{-2\lambda T^{2}},$ for all 
$\lambda \geq \lambda _{1}.$ Also, we keep in mind that by trace theorem 
\[
\left\Vert u\left( \mathbf{x},\pm T\right) \right\Vert _{H^{1}\left( \Omega
\right) },\left\Vert u_{t}\left( \mathbf{x},\pm T\right) \right\Vert
_{L_{2}\left( \Omega \right) }\leq C\left\Vert u\right\Vert
_{H^{k_{n}}\left( Q_{T}^{\pm }\right) },\forall u\in H^{k_{n}}\left(
Q_{T}^{\pm }\right) . 
\]%
Then, taking $\beta \in \left[ 2e^{-\lambda T^{2}},1\right) $ and using (\ref%
{6.8}) and (\ref{6.9}), we obtain%
\begin{equation}
J_{\lambda ,\beta }\left( w_{1}+h\right) -J_{\lambda ,\beta }\left(
w_{1}\right) -J_{\lambda ,\beta }^{\prime }\left( w_{1}\right) \left(
h\right)  \label{6.10}
\end{equation}%
\[
\geq \frac{C_{1}}{\lambda }e^{-2\lambda B^{2}}\mathop{\displaystyle \int}%
\limits_{Q_{T}^{\pm }}\left( h_{t}^{2}+\mathop{\displaystyle \sum }%
\limits_{i,j=1}^{n}h_{x_{i}x_{j}}^{2}\right) \varphi _{\lambda }d\mathbf{x}%
dt+C\lambda e^{-2\lambda B^{2}}\mathop{\displaystyle \int}%
\limits_{Q_{T}^{\pm }}\left[ \left( \nabla h\right) ^{2}+\lambda ^{2}h^{2}%
\right] \varphi _{\lambda }d\mathbf{x}dt 
\]%
\[
+\frac{\beta }{2}\left\Vert h\right\Vert _{H^{k_{n}}\left( Q_{T}^{\pm
}\right) }^{2},\forall h\in \overline{B_{0}\left( 2R\right) },\forall
\lambda \geq \lambda _{1}, 
\]%
also, see (\ref{6.1}). Finally, since $\varphi _{\lambda }\left( x,t\right)
\geq \exp \left( -2\lambda \left( T^{2}-A^{2}\right) \right) $ for $x\in %
\left[ A,B\right] ,t\in \left[ -T,T\right] ,$ then the target estimate (\ref%
{4.201}) follows immediately from (\ref{6.10}). $\ \square $

\subsection{Proof of Theorem 4}

\label{sec:6.2}

Since by (\ref{4.16}) $I_{\lambda ,\beta }\left( W\right) =J_{\lambda ,\beta
}\left( W+G\right) $, $W\in \overline{B_{0}\left( 2R\right) }$ and also
since $W+G\in \overline{B\left( 3R,p_{0},p_{1}\right) },\forall W\in 
\overline{B_{0}\left( 2R\right) },$ then we take in Theorems 2 and 3 $%
\lambda _{1}=\lambda _{1}\left( 3R\right) ,\lambda \left( 3R\right) \geq
\lambda _{1}\left( 3R\right) $ meaning that we replace in (\ref{4.200}) $R$
with $3R$. Denote $w_{1}=W_{1}+G,w_{2}=W_{2}+G.$ Then $w_{1},w_{2}\in 
\overline{B\left( 3R,p_{0},p_{1}\right) }.$ The rest of the proof follows
immediately from Theorems 2 and 3. $\square $

\section{Proof of Theorem 5}

\label{sec:7}

Recall that by (\ref{4.14}) $w^{\ast }-G^{\ast }=W^{\ast }\in B_{0}\left(
2R-\delta \right) .$ Hence, by (\ref{3.13}), (\ref{4.3}), (\ref{4.5}) and (%
\ref{4.7}) $W^{\ast }+G\in B\left( 3R,p_{0},p_{1}\right) .$ We temporally
denote $I_{\lambda ,\beta }\left( W,f_{0}\right) $ the functional $%
I_{\lambda ,\beta }\left( W\right) =J_{\lambda ,\beta }\left( W+G\right) ,$
in which we emphasize the presence of the vector function $\widetilde{f}%
_{0}=\nabla \ln f_{0}$ in the operator $K\left( w\right) .$ We also
temporally denote this operator as $K\left( w,f_{0}\right) =K\left(
W+G,f_{0}\right) $ in (\ref{3.8}) and (\ref{3.15}). Similarly, we also
temporally denote $J_{\lambda ,\beta }\left( W+G,f_{0}\right) :=I_{\lambda
,\beta }\left( W,f_{0}\right) .$ Let 
\begin{equation}
I_{\lambda ,\beta }^{0}\left( W,f_{0}\right) =J_{\lambda ,\beta }^{0}\left(
W+G,f_{0}\right) =e^{-2\lambda B^{2}}\mathop{\displaystyle \int}%
\limits_{Q_{T}^{\pm }}\left( K\left( W+G,f_{0}\right) \right) ^{2}\varphi
_{\lambda }d\mathbf{x}dt.  \label{7.1}
\end{equation}%
By (\ref{3.8}) $K\left( W^{\ast }+G^{\ast },f_{0}^{\ast }\right) =0.$ Hence, 
\begin{equation}
I_{\lambda ,\beta }^{0}\left( W^{\ast },f_{0}^{\ast }\right) =J_{\lambda
,\beta }^{0}\left( W^{\ast }+G^{\ast },f_{0}^{\ast }\right) =0.  \label{7.2}
\end{equation}%
Hence, it follows from (\ref{3.8}), (\ref{4.5}), (\ref{4.50}), (\ref{4.60}),
(\ref{4.8}), (\ref{7.1}) and (\ref{7.2}) that 
\[
I_{\lambda ,\beta }\left( W^{\ast },f_{0}\right) =J_{\lambda ,\beta }\left(
W^{\ast }+G,f_{0}\right) =J_{\lambda ,\beta }\left( W^{\ast }+G^{\ast
}+\left( G-G^{\ast }\right) ,f_{0}^{\ast }+\left( f_{0}-f_{0}^{\ast }\right)
\right) 
\]%
\[
=J_{\lambda ,\beta }^{0}\left( W^{\ast }+G^{\ast },f_{0}^{\ast }\right)
+P_{\lambda ,\delta }+\beta \left\Vert W^{\ast }+G\right\Vert
_{H^{k_{n}}\left( Q_{T}^{\pm }\right) }=P_{\lambda ,\delta }+\beta
\left\Vert W^{\ast }+G\right\Vert _{H^{k_{n}}\left( Q_{T}^{\pm }\right)
}^{2}, 
\]%
where $\left\vert P_{\lambda ,\delta }\right\vert \leq C_{1}\delta ^{2}.$
Thus, 
\begin{equation}
I_{\lambda ,\beta }\left( W^{\ast },f_{0}\right) \leq C_{1}\delta ^{2}+\beta
\left\Vert W^{\ast }+G\right\Vert _{H^{k_{n}}\left( Q_{T}^{\pm }\right)
}^{2}.  \label{7.3}
\end{equation}

By (\ref{4.5}) and (\ref{4.6}) 
\[
\left\Vert W^{\ast }+G\right\Vert _{H^{k_{n}}\left( Q_{T}^{\pm }\right)
}=\left\Vert \left( W^{\ast }+G^{\ast }\right) +\left( G-G^{\ast }\right)
\right\Vert _{H^{k_{n}}\left( Q_{T}^{\pm }\right) }\leq \left\Vert w^{\ast
}\right\Vert _{H^{k_{n}}\left( Q_{T}^{\pm }\right) }+\delta <R. 
\]%
Hence, 
\begin{equation}
W^{\ast }+G\in B\left( R,p_{0},p_{1}\right) .  \label{7.4}
\end{equation}%
Let $W_{\min ,\lambda \left( 3R\right) ,\beta \left( 3R\right) }\in 
\overline{B_{0}\left( 2R\right) }$ be the minimizer of the functional $%
I_{\lambda ,\beta }\left( W,f_{0}\right) ,$ the existence and uniqueness of
which on the set $\overline{B_{0}\left( 2R\right) }$ is guaranteed by
Theorem 3. We will choose the dependencies on $\delta $ of parameters $%
\lambda $ and $\beta $ later in this proof. Thus, we can apply (\ref{4.17})
now as%
\[
I_{\lambda \left( 3R\right) ,\beta \left( 3R\right) }\left( W^{\ast
},f_{0}\right) -I_{\lambda \left( 3R\right) ,\beta \left( 3R\right) }\left(
W_{\min ,\lambda \left( 3R\right) ,\beta \left( 3R\right) },f_{0}\right) 
\]%
\[
-I_{\lambda \left( 3R\right) ,\beta \left( 3R\right) }^{\prime }\left(
W_{\min ,\lambda \left( 3R\right) ,\beta \left( 3R\right) },f_{0}\right)
\left( W^{\ast }-W_{\min ,\lambda \left( 3R\right) ,\beta \left( 3R\right)
}\right) 
\]%
\begin{equation}
\geq C_{1}\exp \left( -3\lambda \left( 3R\right) \left( \gamma
^{2}T^{2}+B^{2}-A^{2}\right) \right) \left\Vert W^{\ast }-W_{\min ,\lambda
\left( 3R\right) ,\beta \left( 3R\right) }\right\Vert _{H^{2,1}\left(
Q_{T\gamma }^{\pm }\right) }^{2}  \label{7.5}
\end{equation}%
\[
+\frac{\beta \left( 3R\right) }{2}\left\Vert W^{\ast }-W_{\min ,\lambda
\left( 3R\right) ,\beta \left( 3R\right) }\right\Vert _{H^{k_{n}}\left(
Q_{T\gamma }^{\pm }\right) }^{2}. 
\]%
By (\ref{4.170}) $-I_{\lambda \left( 3R\right) ,\beta \left( 3R\right)
}^{\prime }\left( W_{\min ,\lambda \left( 3R\right) ,\beta \left( 3R\right)
},f_{0}\right) \left( W^{\ast }-W_{\min ,\lambda \left( 3R\right) ,\beta
\left( 3R\right) }\right) \leq 0.$ Recall that $\eta _{1}=\gamma
^{2}T^{2}+B^{2}-A^{2}.$ Hence, (\ref{7.3})-(\ref{7.5}) imply that%
\begin{equation}
\left\Vert W^{\ast }-W_{\min ,\lambda \left( 3R\right) ,\beta \left(
3R\right) }\right\Vert _{H^{2,1}\left( Q_{T\gamma }^{\pm }\right) }^{2}\leq
C_{1}\left( \delta ^{2}+\beta \right) \exp \left( 3\lambda \left( 3R\right)
\eta _{1}\right) .  \label{7.7}
\end{equation}

We now specify dependencies of $\lambda $ and $\beta $ on the noise level $%
\delta .$ Choose $\lambda =\lambda \left( \delta ,3R\right) $ such that 
\[
\exp \left( 3\lambda \left( \delta ,3R\right) \eta _{1}\right) =\exp \left(
3\lambda \left( \delta ,3R\right) \left( \gamma ^{2}T^{2}+B^{2}-A^{2}\right)
\right) =\frac{1}{\delta }. 
\]%
Then $\lambda \left( \delta ,3R\right) =\ln \left( \delta ^{-1/\left( 3\eta
_{1}\right) }\right) $ and 
\begin{equation}
\delta ^{2}\exp \left( 3\lambda \left( \delta ,3R\right) \eta _{1}\right)
=\delta .  \label{7.8}
\end{equation}%
Since $\delta \in \left( 0,\delta _{0}\right) $ and since $\ln \left( \delta
_{0}^{-1/\left( 3\eta _{1}\right) }\right) \geq \lambda _{1}\left( 3R\right)
,$ then $\lambda \left( \delta ,3R\right) \geq \lambda _{1}\left( 3R\right)
. $ Next, by Theorem 2, we can take $\beta =2e^{-\lambda \left( \delta
,3R\right) T^{2}}.$ Hence, in (\ref{7.7}) 
\begin{equation}
\beta \exp \left( 3\lambda \left( 3R\right) \eta _{1}\right) =\delta ^{\eta
_{2}/\eta _{1}},\eta _{2}=\left( 1-3\gamma ^{2}\right) T^{2}-3\left(
B^{2}-A^{2}\right) >0.  \label{7.9}
\end{equation}%
Recalling that $2\rho =\min \left( 1,\eta _{2}/\eta _{1}\right) $ and using (%
\ref{7.7})-(\ref{7.9}), we obtain%
\[
\left\Vert W^{\ast }-W_{\min ,\lambda \left( 3R\right) ,\beta \left(
3R\right) }\right\Vert _{H^{2,1}\left( Q_{\gamma T}^{\pm }\right) }\leq
C_{1}\delta ^{\rho }. 
\]%
Hence,%
\[
\left\Vert w^{\ast }-w_{\min ,\lambda \left( \delta ,3R\right) ,\beta \left(
\delta ,3R\right) }\right\Vert _{H^{2,1}\left( Q_{\gamma T}^{\pm }\right)
}\leq \left\Vert W^{\ast }-W_{\min ,\lambda \left( \delta ,3R\right) ,\beta
\left( \delta ,3R\right) }\right\Vert _{H^{2,1}\left( Q_{\gamma T}^{\pm
}\right) } 
\]%
\begin{equation}
+\left\Vert G^{\ast }-G\right\Vert _{H^{2,1}\left( Q_{\gamma T}^{\pm
}\right) }\leq C_{2}\delta ^{\rho }+\delta \leq \left( C_{2}+1\right) \delta
^{\rho },  \label{7.10}
\end{equation}%
which proves (\ref{4.18}). Finally, since (\ref{4.18}) holds, then (\ref%
{4.19}) follows immediately from (\ref{3.16}) and the rest of the discussion
in the last paragraph of section 3. \ $\square $

\section{Proof of Theorem 7}

\label{sec:8}

Recall that Theorem 6 is valid: see item 2 in Remarks 3 (section 4). By the
triangle inequality, (\ref{4.5}), (\ref{4.18}) and (\ref{4.24})%
\[
\left\Vert w^{\ast }-w_{n}\right\Vert _{H^{2,1}\left( Q_{\gamma T}^{\pm
}\right) }=\left\Vert w^{\ast }-w_{\min ,\lambda \left( \delta ,3R\right)
,\beta \left( \delta ,3R\right) }+\left( w_{\min ,\lambda \left( \delta
,3R\right) ,\beta \left( \delta ,3R\right) }-w_{n}\right) \right\Vert
_{H^{2,1}\left( Q_{\gamma T}^{\pm }\right) } 
\]%
\[
\leq C_{2}\delta ^{\rho }+\left\Vert w_{\min ,\lambda \left( \delta
,3R\right) ,\beta \left( \delta ,3R\right) }-w_{n}\right\Vert
_{H^{2,1}\left( Q_{\gamma T}^{\pm }\right) }\leq C_{2}\delta ^{\rho
}+\left\Vert w_{\min ,\lambda \left( \delta ,3R\right) ,\beta \left( \delta
,3R\right) }-w_{n}\right\Vert _{H^{2,1}\left( Q_{T}^{\pm }\right) } 
\]%
\[
=C_{2}\delta ^{\rho }+\left\Vert W_{\min ,\lambda \left( 3R\right) ,\beta
}-W_{n}\right\Vert _{H^{k_{n}}\left( Q_{T}^{\pm }\right) } 
\]%
\[
\leq C_{2}\delta ^{\rho }+\theta ^{n}\left\Vert W_{\min ,\lambda \left(
3R\right) ,\beta }-W_{0}\right\Vert _{H^{k_{n}}\left( Q_{T}^{\pm }\right)
}=C_{2}\delta ^{\rho }+\theta ^{n}\left\Vert w_{\min ,\lambda \left(
3R\right) ,\beta }-w_{0}\right\Vert _{H^{k_{n}}\left( Q_{T}^{\pm }\right) }, 
\]%
which proves (\ref{4.25}). Estimate (\ref{4.26}) follows immediately from (%
\ref{4.25}) and the discussion in the last paragraph of section 3. \ $%
\square $

\section{Numerical Testing}

\label{sec:9}

In the following tests, we set the domain $\Omega =\left( 1,2\right) \times
\left( 1,2\right) $ and also%
\[
Lu=\Delta u-c\left( \mathbf{x}\right) u. 
\]
To solve the inverse problem, we should first computationally simulate the
data (\ref{2.7}), (\ref{2.8}) via the numerical solution of the forward
problem (\ref{2.4}). To solve problem (\ref{2.4}), computationally, we have
used the standard finite difference method. The spatial mesh size is $%
1/640\times 1/640$ while the temporal one $T/512$. For the forward problem,
we use the implicit scheme to compute the data needed for the inverse
problem.

In computations of the inverse problem, the spatial mesh size is $1/16\times
1/16$ and the temporal one $T/16$. When minimizing the functional $%
J_{\lambda ,\beta }\left( w\right) $ in the discrete sense, we formulate the
right hand side of (\ref{3.15}) via finite differences and minimize with
respect to the values of the function $w$ at grid points. To minimize the
discretized functional, 
we use Matlab's built-in function \textbf{fminunc} with its option of 
\textit{quasi-newton algorithm}. 
This procedure calculates the gradient $\nabla J_{\lambda ,\beta }\left(
w\right) $ automatically and iterations stop when the condition $\left\vert
\nabla J_{\lambda ,\beta }\left( w\right) \right\vert <1\times 10^{-2}$
holds. Note that even though our theory requires the application of the
gradient projection method, we have established numerically that we can
avoid the use of the projection operator $P_{B}$ and to use just the
conjugate gradient method. In fact, the use of the operator $P_{B}$ would
complicate the matter. The same observation took place in all our works on
the convexification, which contain numerical studies \cite%
{Khoa,convIPnew,EIT,timedomain,KlibKol3}. Also, we have minimized the
functional $J_{\lambda ,\beta }\left( w\right) $ rather than $I_{\lambda
,\beta }\left( W\right) $ and it worked quite well.

As to (\ref{3.16}), we have numerically discovered that rather than taking
an average over $t\in \left[ -\gamma T,\gamma T\right] ,$ better to use (\ref%
{3.4}) at $\left\{ t=t_{0}\right\} .$ In numerical tests below, we took 
\begin{equation}
\lambda =1,\text{ }k_{n}=3,\text{ }\beta =0.01.  \label{9.1}
\end{equation}%
In the process of the minimization of the functional $J_{\lambda
,\beta }\left( w\right) ,$ the starting point of iterations is always chosen to be 
the null function of value zero everywhere.

In the following three tests, we show the results of the recovery of the
coefficients $c\left( \mathbf{x}\right) $ with sophisticated structures. We
choose the tested coefficients $c\left( \mathbf{x}\right) $ having the
shapes of the letters `$A$' and `$\Omega $'. We measure $g_{1}(x_{1},x_{2},t)
$ on $16\times 32$ detectors uniformly distributed on the rectangle $\Gamma
_{T}^{\pm }$ and `measure' the function $f_{0}(x_{1},x_{2},t_{0})$ on $%
16\times 16$ detectors uniformly distributed on the square $(1,2)\times
(1,2)\times \left\{ t=t_{0}\right\} $. As initial and Dirichlet boundary
conditions for the data simulations in (\ref{2.5}), (\ref{2.6}), we took 
\[
u\left( \mathbf{x},-T\right) =1+\sin (\pi (x_{1}-1))\sin (\pi (x_{2}-1))%
\text{ and }u\mid _{S_{T}^{\pm }}=1. 
\]%
We allow in our tests the function $c\left( \mathbf{x}\right) $\ to be both
positive and negative. Indeed, we have imposed the positivity condition (\ref%
{2.30}) only to ensure that the function $u\left( \mathbf{x},t\right) \neq 0$%
\ in $\overline{Q_{T}^{\pm }}.$\ However, we have not observed any zeros of
this function in our numerical studies.

\smallskip

\begin{figure}[tbp]
\begin{center}
\begin{tabular}{cc}
\includegraphics[width=4cm]{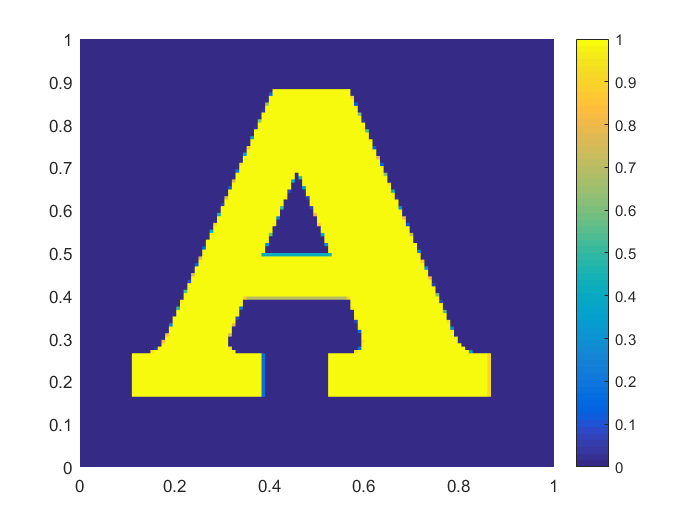} & %
\includegraphics[width=4cm]{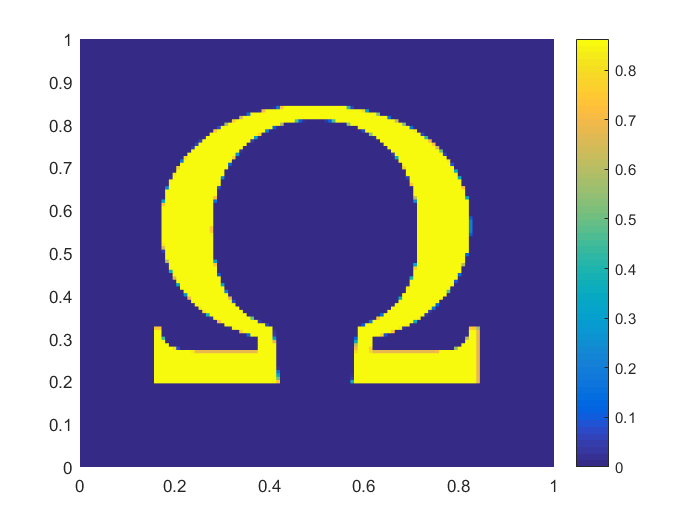} \\ 
(a) $c\left( \mathbf{x}\right)$ with the shape of the letter '$A$' & (d) $%
c\left( \mathbf{x}\right)$ with the shape of the letter '$\Omega$' \\ 
\includegraphics[width=4cm]{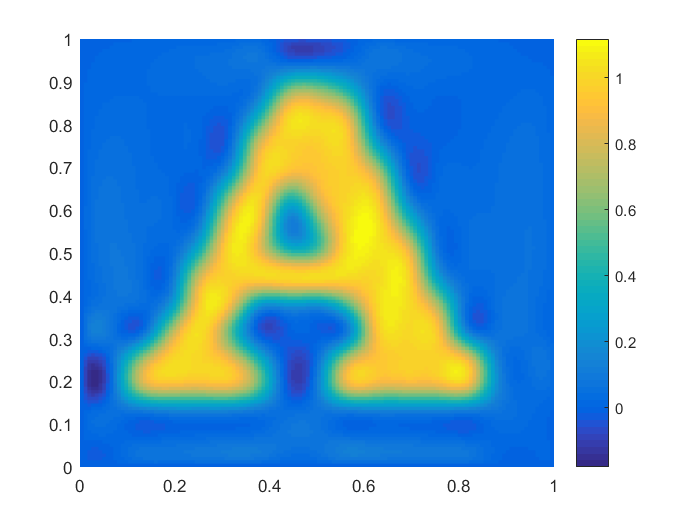} & %
\includegraphics[width=4cm]{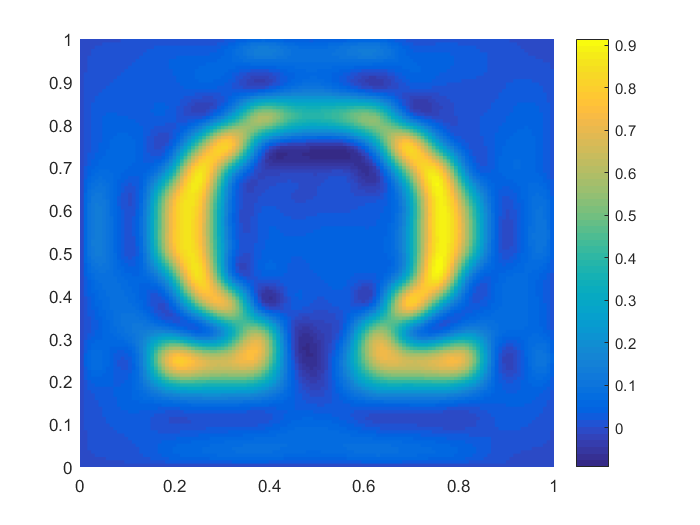} \\ 
(b) Recovered $c\left( \mathbf{x}\right)$ for $T=1$ & (e) Recovered $c\left( 
\mathbf{x}\right)$ for $T=1$ \\ 
\includegraphics[width=4cm]{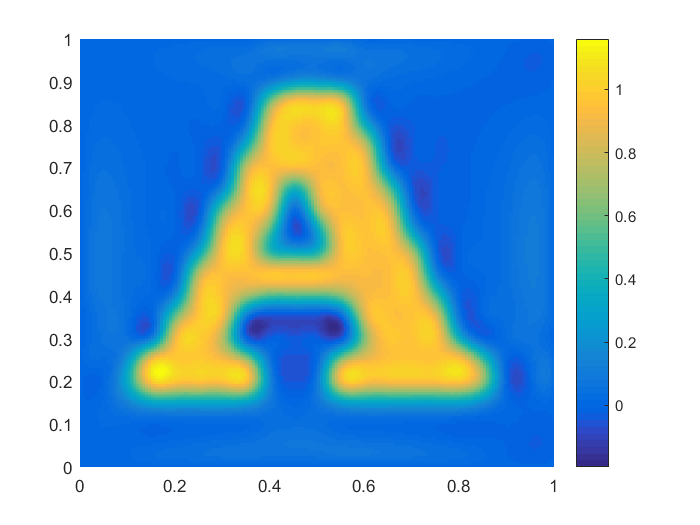} & %
\includegraphics[width=4cm]{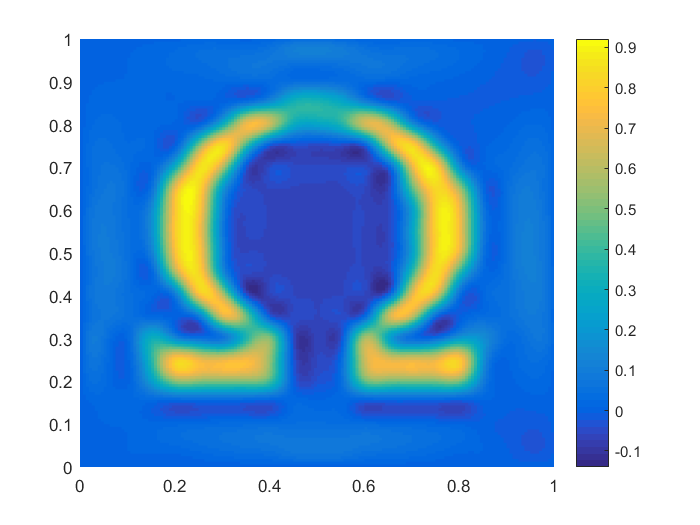} \\ 
(c) Recovered $c\left( \mathbf{x}\right)$ for $T=0.1$ & (f) Recovered $%
c\left( \mathbf{x}\right)$ for $T=0.1$%
\end{tabular}%
\end{center}
\caption{\emph{Results of Test 1. Here $t_{0}=0$ in (\protect\ref{2.8}). (a)
The coefficient $c\left( \mathbf{x}\right)$ with the shape of the letter
'A'. (d) The coefficient $c\left( \mathbf{x}\right)$ with the shape of the
letter '$\Omega$'. (b) and (c) are the recovered $c\left( \mathbf{x}\right)$
for $T=1$ and $T=0.1$ respectively for coefficient with the shape of the
letter 'A'. (e) and (f) are the recovered $c\left( \mathbf{x}\right)$ for $%
T=1$ and $T=0.1$ respectively for coefficient with the shape of the letter '$%
\Omega$'.}}
\label{example1}
\end{figure}

\textbf{Test 1}. First, we test the reconstruction by our method of the
coefficients $c\left( \mathbf{x}\right) $ with the shapes of letters `$A$'
and `$\Omega $'. In this test, we measure the data at time $\left\{
t_{0}=0\right\} $ for the cases $T=1$ and $T=0.1$. The numerical results are
shown in Figure \ref{example1}.

\smallskip

\begin{figure}[tbp]
\begin{center}
\begin{tabular}{cc}
\includegraphics[width=4.8cm]{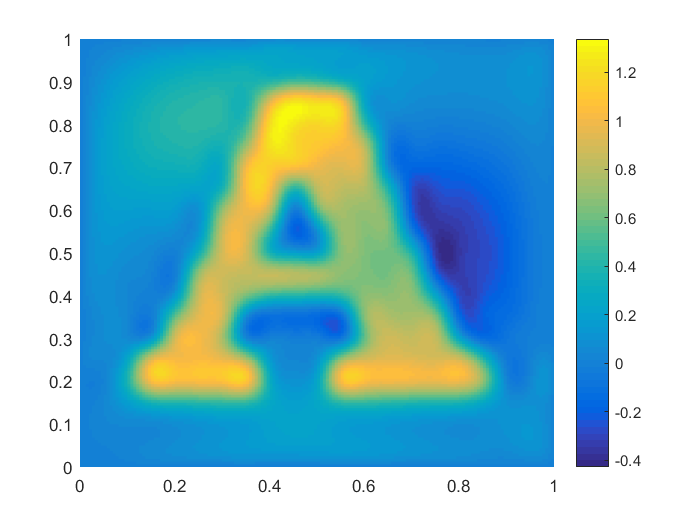} & %
\includegraphics[width=4.8cm]{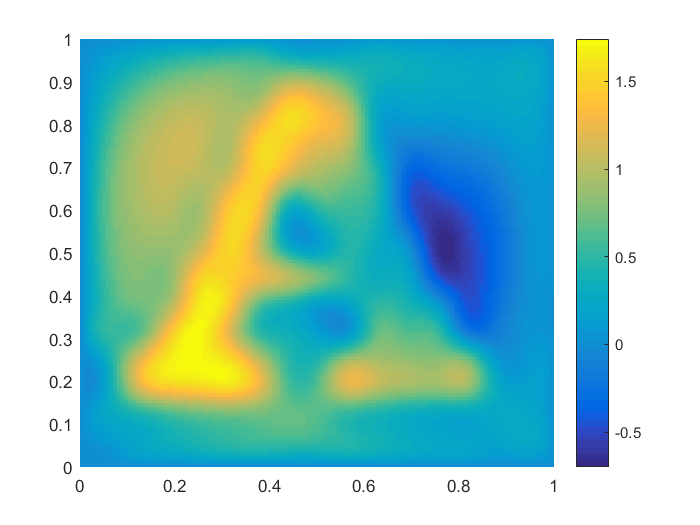} \\ 
(a)Recovered $c\left( \mathbf{x}\right)$ for $\epsilon=0.02$ & (b) Recovered 
$c\left( \mathbf{x}\right)$ for $\epsilon=0.01$ \\ 
\includegraphics[width=4.8cm]{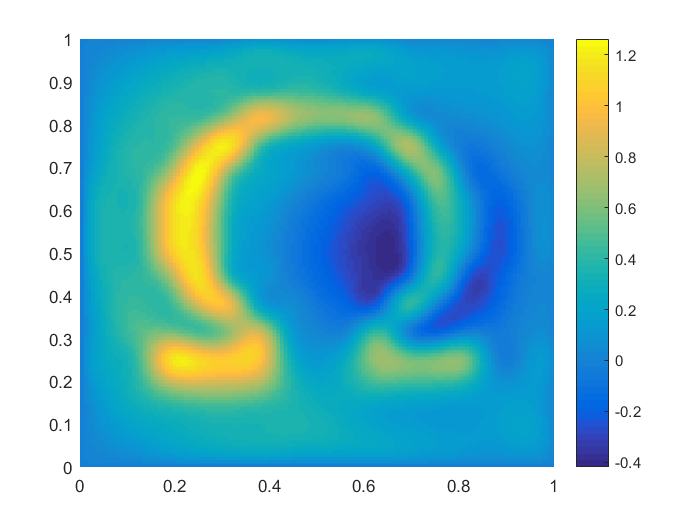} & %
\includegraphics[width=4.8cm]{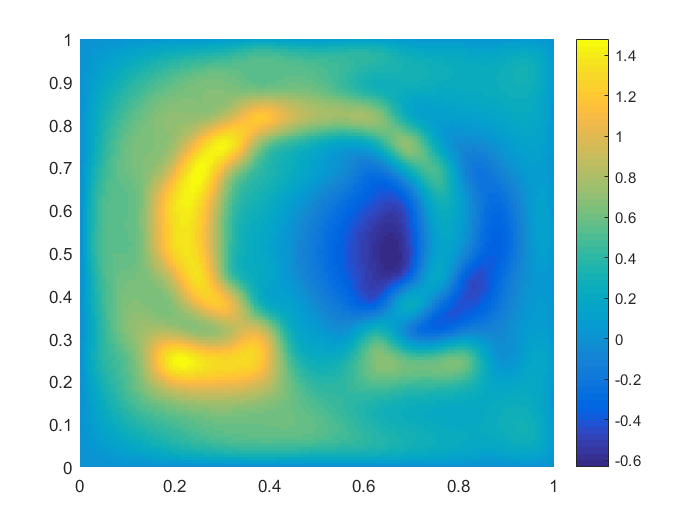} \\ 
(c) Recovered $c\left( \mathbf{x}\right)$ for $\epsilon=0.02$ & (d)
Recovered $c\left( \mathbf{x}\right)$ for $\epsilon=0.01$%
\end{tabular}%
\end{center}
\caption{\emph{Results of Test 2. Here $T=0.1$ and in (\protect\ref{2.8}) $%
t_{0}=0,$. We show the results in the case when the data are measured at a
time $t_{0},$ which is close to the initial time moment $\left\{
t=-T=-0.1\right\}$. (a) and (b) are the recovered $c\left( \mathbf{x}\right)$
for $\protect\epsilon=0.02$ and $\protect\epsilon=0.01$ respectively for
coefficient with the shape of the letter 'A'. (c) and (d) are the recovered $%
c\left( \mathbf{x}\right)$ for $\protect\epsilon=0.02$ and $\protect\epsilon%
=0.01$ respectively for coefficient with the shape of the letter '$\Omega$'.
Comparison with Figure 1 shows that the quality of images is better if $%
t_{0} $ is not too close to the initial time moment $\left\{ t=-T\right\}$.}}
\label{example2}
\end{figure}

\textbf{Test 2}. In this test, we set $T=0.1$. We show the results in the
case when the data are measured at a time $\left\{ t_{0}\right\} $ which is
close to the initial time $\left\{ t=-T=-0.1\right\} $. We take $%
t_{0}=-T+\epsilon $ with $\epsilon =0.02$ and $\epsilon =0.01$. We test the
reconstruction by our method of the coefficients $c\left( \mathbf{x}\right) $
with the shapes of the letters `$A$' and `$\Omega $'. The numerical results
are shown in Figure \ref{example2}. In this test, we demonstrate the results
when one measures the data at some time close to the initial time. It is
numerically shown that the closer $t_{0}$ is to the initial time $t=-T$, the
worse the result is.

\smallskip

\begin{figure}[tbp]
\begin{center}
\begin{tabular}{cc}
\includegraphics[width=4.8cm]{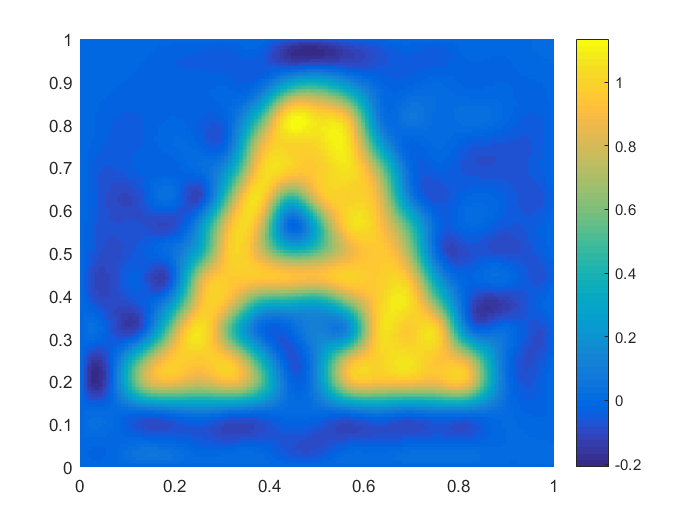}
& %
\includegraphics[width=4.8cm]{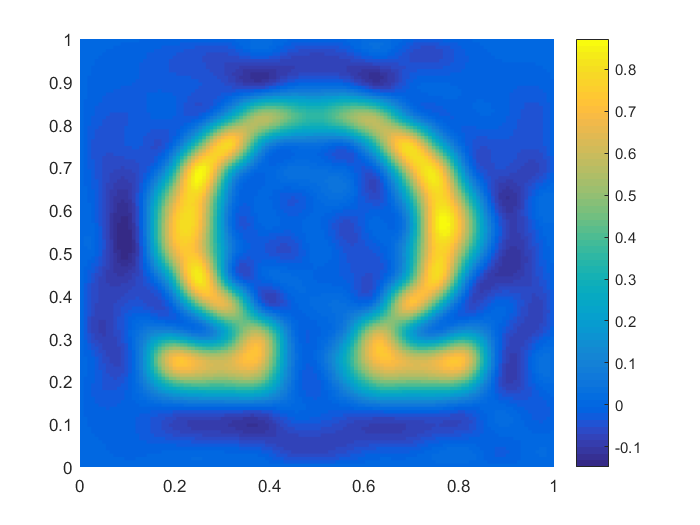}
\\ 
(a)Recovered $c\left( \mathbf{x}\right)$ with the shape 'A' & (b) Recovered $%
c\left( \mathbf{x}\right)$ with the shape '$\Omega $'%
\end{tabular}%
\end{center}
\caption{\emph{Results of Test 3. Here $T=1$ and $t_{0}=0$ in (\protect\ref%
{2.8}). (a) and (b) are the recovered coefficients $c\left( \mathbf{x}%
\right) $ with the shapes of the letters '$A$' and '$\Omega$' respectively.
The measured data contain $5\%$ relative random noise.}}
\label{example3}
\end{figure}

\textbf{Test 3}. We now want to see how the random noise in the data
influences our reconstruction. We add $5\%$ relative random noise to each
detector on $\Gamma _{T}^{\pm }$ as well as on $(1,2)\times (1,2)\times
\left\{ t=0\right\} $, i.e. we work now with the noisy data, 
\begin{equation}
u_{x}^{noise}\mid _{\Gamma _{T}^{\pm }}=g_{1}\left( \mathbf{x},t\right)
+\sigma \xi _{\mathbf{x},t}g_{1}\left( \mathbf{x},t\right) ,  \label{9.2}
\end{equation}%
\begin{equation}
u^{noise}\left( \mathbf{x},t_{0}\right) =f_{0}\left( \mathbf{x}\right)
+\sigma \xi _{\mathbf{x}}f_{0}\left( \mathbf{x}\right) .  \label{9.3}
\end{equation}%
Here $\sigma =5\%$ is the noise level, $\xi _{\mathbf{x},t}$ and $\xi _{%
\mathbf{x}}$ are independent normally distributed random variables. To
preprocess the noisy data, we use the thin plate spline smoother developed
in \cite{ZhangChen18}. The algorithm proposed in \cite{ZhangChen18} provides
a good approximation to the true function without knowing neither the noise
level nor any other a priori information of the true function to be
approximated. Then the cubic B-splines are employed to approximate the first
and second order derivatives of the noisy data. In this test, we `measure' $%
g_{1}(x_{1},x_{2},t)$ on $16\times 32$ detectors uniformly distributed on
the plane $\Gamma _{T}^{\pm }$ and also `measure' $f_{0}(x_{1},x_{2},t_{0})$
on $160\times 160$ detectors uniformly distributed on the plane $(1,2)\times
(1,2)\times \left\{ t=0\right\} $. We now set $T=1$, in (\ref{2.8}) $%
t_{0}=0, $and the noise is added to the data as in (\ref{9.2}), (\ref{9.3}).
We test the reconstruction by our method of the coefficients with the shape
of the letters `$A$' and `$\Omega $'. The numerical results are shown in
Figure \ref{example3}. We see that our method works still very well in the
mild noisy case.

\end{document}